\documentclass{iopart}

\usepackage{t1enc}
\usepackage{graphicx}
\usepackage[english]{babel}
\usepackage[T1]{fontenc}
\usepackage[latin1]{inputenc}
\usepackage{url}
\usepackage{dsfont}
\usepackage{graphicx}
\usepackage{graphics}
\usepackage{subfig}
\usepackage{iopams} 
\usepackage{setstack}
\usepackage{amstext}
\usepackage{amsopn}
\usepackage{amsthm}
\usepackage{enumitem}
\usepackage{color}

\newtheorem{algo}{Algorithm}

\DeclareMathOperator*{\mydef}{\mathrel{\mathop:}=}
\DeclareMathOperator*{\mydefswitched}{=\mathrel{\mathop:}}
\newcommand{\R}{\mathbb{R}}
\newcommand{\J}{\mathcal{J}}
\newcommand{\norm}[1]{\Vert #1 \Vert}
\newcommand{\sqnorm}[1]{\Vert #1 \Vert_2^2}

\newcommand{\abs}[1]{\vert #1 \vert}

\newcommand{\argmin}{\text{argmin}}
\newcommand{\argmaxsub}[1]{\underset{{ #1 }}{{\rm argmax}}}

\newcommand{\Exp}{\mathbb{E}}

\newcommand{\uMAP}{\hat{u}_{\text{\tiny MAP}}}
\newcommand{\uCM}{\hat{u}_{\text{\tiny CM}}}

\newcommand{\post}{p_{post}(u|f)}
\newcommand{\prior}{p_{prior}(u)}
\newcommand{\like}{p_{like}(f|u)}
\newcommand{\Indicator}{\mathds{1}}

\newcommand{\termabb}[2]{\emph{#1} (\emph{#2})}
\newcommand{\proptoab}[1]{\stackrel{ #1 }{\propto}}
\newcommand{\uInf}{u^{\infty}}
\newcommand{\Radon}{\mathcal{R}}

\begin{document}

\title[]{Fast Gibbs sampling for high-dimensional Bayesian inversion.}

\author{Felix Lucka}

\address{Centre for Medical Image Computing, University College London, WC1E 6BT London, UK}

\ead{f.lucka@ucl.ac.uk}

\begin{abstract}
Solving ill-posed inverse problems by Bayesian inference has recently attracted considerable attention. Compared to deterministic approaches, the probabilistic representation of the solution by the posterior distribution can be exploited to explore and quantify its uncertainties. In applications where the inverse solution is subject to further analysis procedures, this can be a significant advantage. Alongside theoretical progress, various new computational techniques allow to sample very high dimensional posterior distributions: In \cite{Lu12}, a Markov chain Monte Carlo (MCMC) posterior sampler was developed for linear inverse problems with $\ell_1$-type priors. In this article, we extend this single component Gibbs-type sampler to a wide range of priors used in Bayesian inversion, such as general $\ell_p^q$ priors with additional hard constraints. Besides a fast computation of the conditional, single component densities in an explicit, parameterized form, a fast, robust and exact sampling from these one-dimensional densities is key to obtain an efficient algorithm. We demonstrate that a generalization of slice sampling can utilize their specific structure for this task and illustrate the performance of the resulting slice-within-Gibbs samplers by different computed examples. These new samplers allow us to perform sample-based Bayesian inference in high-dimensional scenarios with certain priors for the first time, including the inversion of computed tomography (CT) data with the popular isotropic total variation (TV) prior.
\end{abstract}

\ams{65J22,62F15,65C05,65C60}
\submitto{\IP}
\maketitle


\section{Introduction} \label{sec:Intro}

\subsection{Bayesian Inversion} \label{subsec:setting}

\begin{figure}[tb]
   \centering
\subfloat[][$\prior \propto \exp \left( - \lambda \sqnorm{u} \right)$ \label{subfig:L2LiPrPo}]{\includegraphics[width=0.325\textwidth]{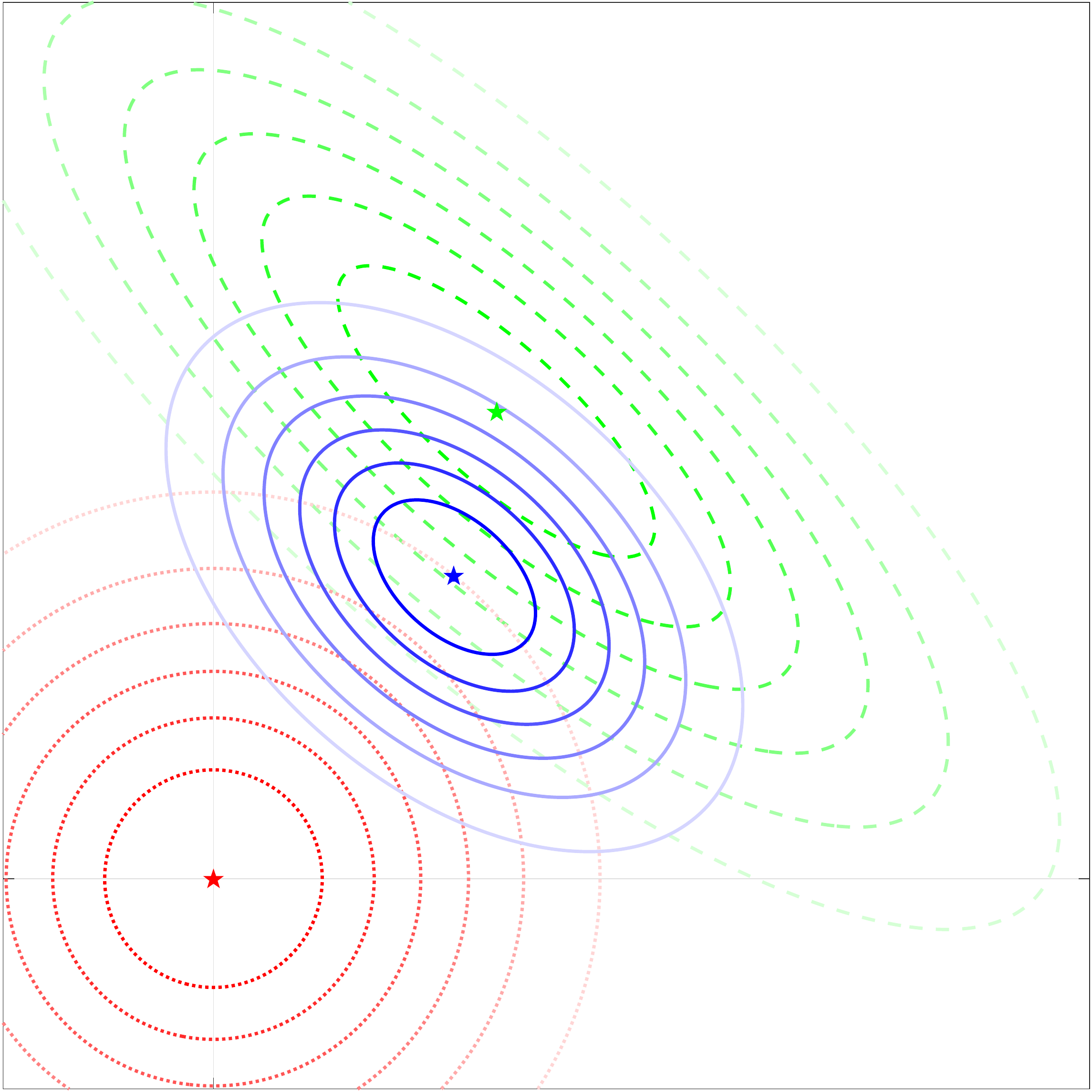}}
\hskip 0.005\textwidth
\subfloat[][$\prior \propto \exp \left( - \lambda \norm{u}_1 \right)$ \label{subfig:L1LiPrPo}]{\includegraphics[width=0.325\textwidth]{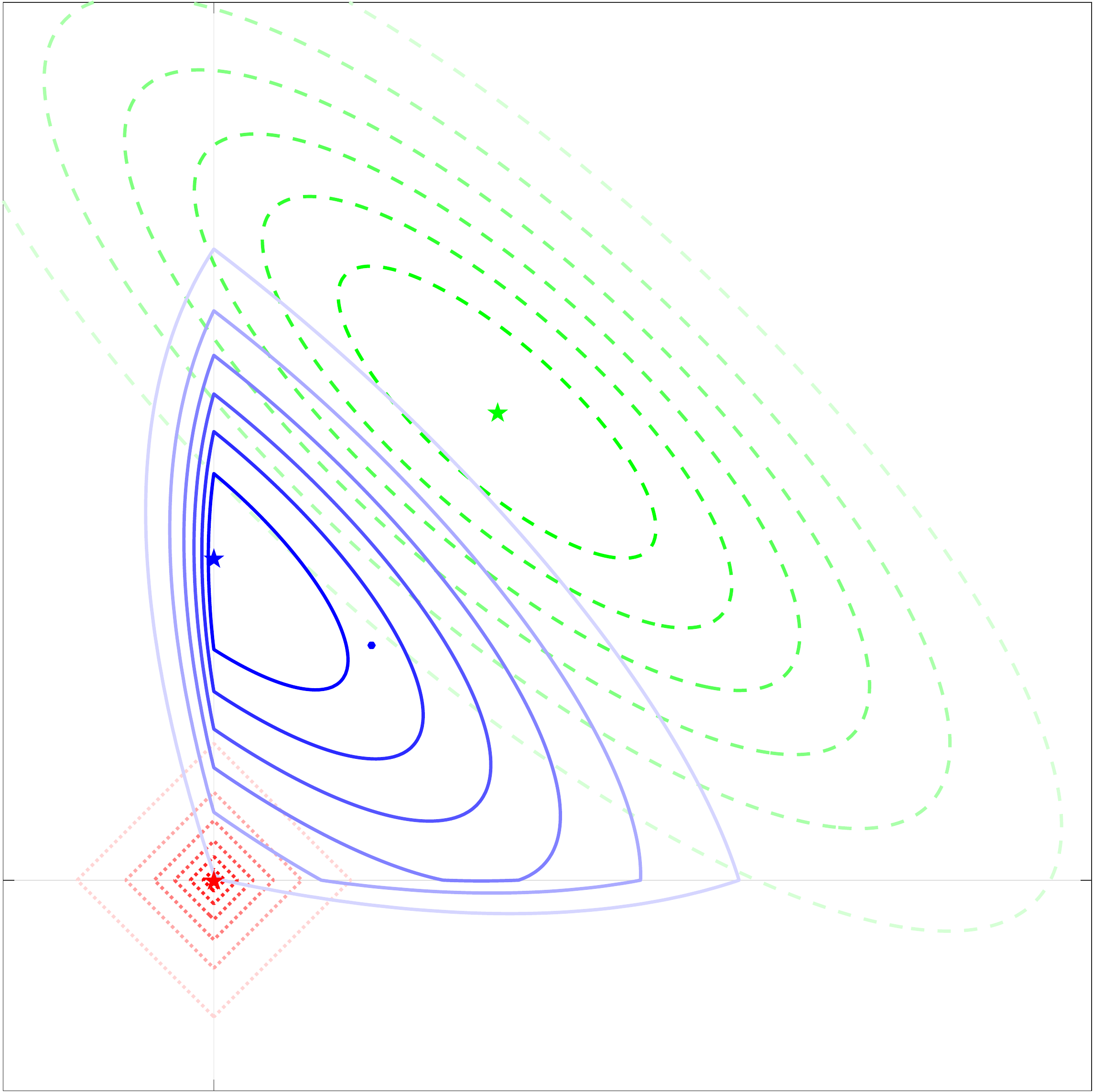}}
\hskip 0.005\textwidth
\subfloat[][SC Gibbs sampler\label{subfig:Gibbs}]{\includegraphics[width=0.3225\textwidth]{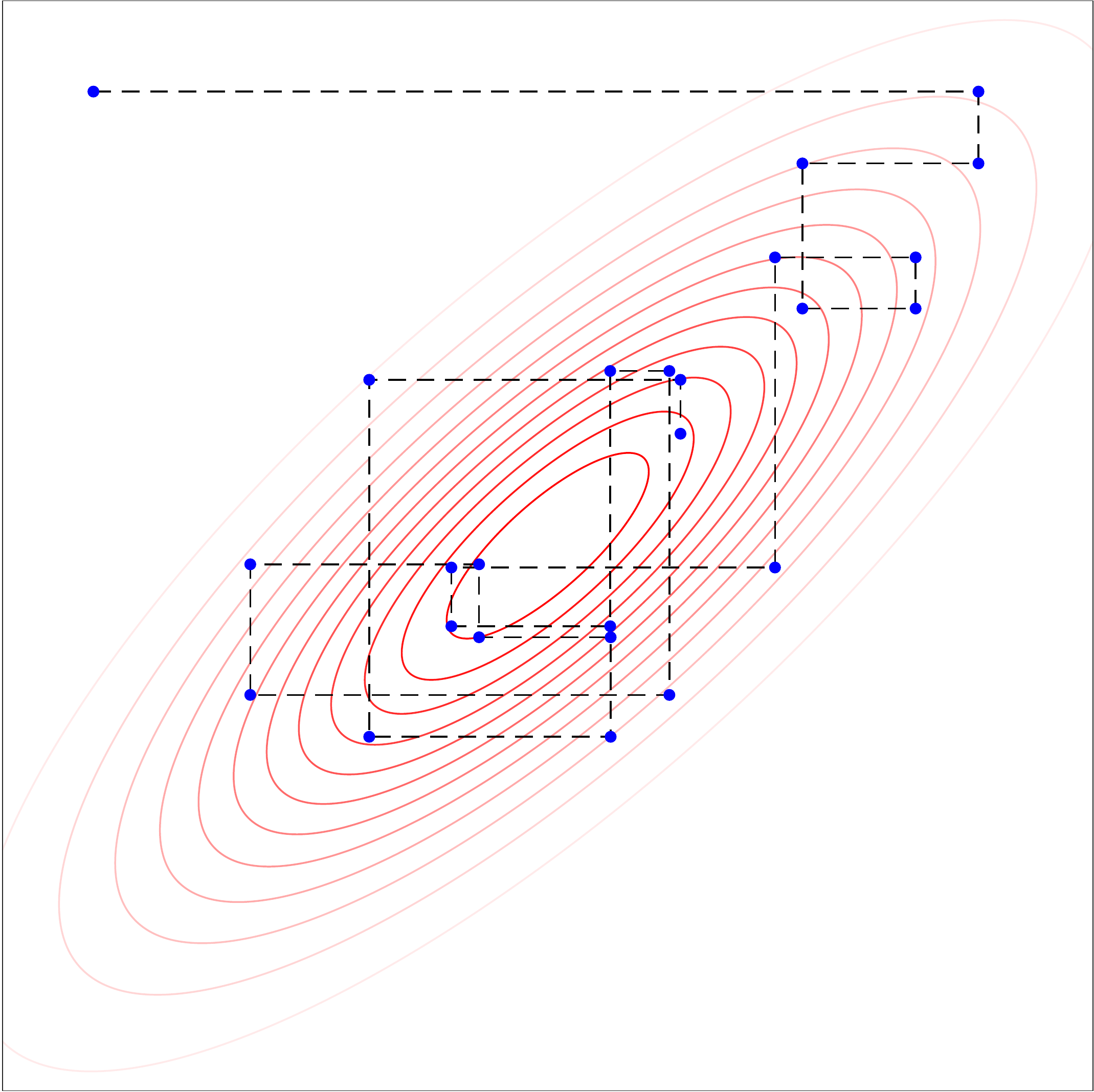}}
\caption{\protect\subref{subfig:L2LiPrPo}-\protect\subref{subfig:L1LiPrPo} Illustration of Bayesian inversion with different priors. Depicted are the level sets of likelihood (green, dashed), prior (red, dotted) and resulting posterior distribution (blue, solid). Maximal and expected value of the corresponding distributions are depicted by a star and a bullet, respectively. \protect\subref{subfig:Gibbs} Example of an MCMC chain generated by an SC Gibbs sampler (blue bullets connected by dashed lines)  to sample a bivariate target density (level sets shown as red solid lines).}
   \label{fig:LiPrPo1}
\end{figure}

We consider the task of inferring information about an unknown quantity from indirect, noisy measurements and assume that a reasonable mathematical model is given by a linear, ill-posed operator equation including additive noise terms. The following discrete forward model is used to carry out the computational inference:
\begin{equation}
  f = A \, u + \varepsilon. \label{eq:FwdEq}
\end{equation}   
Here, $f \in \R^m$ represents the measured data, $u \in \R^n$ a suitable discretization of the unknown quantity we wish to reconstruct and $A \in \R^{m \times n}$ a corresponding discretization of the continuous forward operator. We assume that the statistics of the additive noise can be well-approximated by a Gaussian distribution $\mathcal{N}(\mu,\Sigma)$ and that $f$ and $A$ are already centered and decorrelated with respect to $\mu$ and $\Sigma$, i.e., $f = \Sigma^{-1/2}( \tilde{f}-\mu )$ and $A = \Sigma^{-1/2} \tilde{A}$, where $\tilde{f}$ and $\tilde{A}$ denote the original data and forward operator, respectively. This leads to $\varepsilon \sim \mathcal{N}(0,I_m)$ and the following \emph{likelihood} distribution,
\begin{equation}
  \like \propto \exp \left( -\case{1}{2} \sqnorm{ f - A \, u }\right), \label{eq:Likelihood}
\end{equation}
which is a probabilistic model of the measured data $f$ given the unknown solution $u$. In typical inverse problems, solving \eref{eq:FwdEq} for $u$ is ill-posed. As a consequence, the information that \eref{eq:Likelihood} contains about $u$ is insufficient to carry out robust inference and we need to amend it by \emph{a-priori} information, encoded in a \emph{prior} distribution $\prior$. Then, the total information on $u$ we have \emph{after} performing the measurement is encoded by the conditional distribution of $u$ given $f$, the so-called \emph{a-posteriori} distribution. It can be computed by \emph{Bayes' rule}: 
\begin{equation}
 \post = \frac{\like \prior}{p(f)} \label{eq:BayesRule}
\end{equation}
Figure \ref{fig:LiPrPo1} illustrates the inference process. Originating from statistical physics, \emph{Gibbs distributions} are commonly used prior models:  
\begin{equation}
\prior \propto \exp \left(- \lambda \J(u) \right)  \label{eq:GibbsPrior}
\end{equation} 
The functional $\J(u)$ measures an \emph{energy} of $u$. The use of Gibbs priors leads to
\begin{equation}
\post \propto \exp \left( -\case{1}{2} \sqnorm{ f - A \, u } - \lambda \J(u) \right).
\end{equation}
For a general introduction to \emph{Bayesian inversion} we refer to \cite{KaSo05,St10,Lu14} and references therein. 
The recent attention on this particular inversion approach is best reflected by the recent special issue of \IP \cite{CaKaSo14}, which also provides a good overview over current developments and trends.

\subsection{Sample-based Inference} \label{subsec:SampleBasedInference}

While the posterior $\post$ represents our complete knowledge about $u$, \emph{Bayesian estimation} tries to extract the information of interest from it. Classical examples thereof include the \termabb{maximum a-posteriori estimate}{MAP} and the and the \termabb{conditional mean estimate}{CM},
\begin{equation}
\fl \quad \uMAP \mydef \argmaxsub{u \in \R^n} \left\lbrace   \; p_{post}(u|f) \right\rbrace, \qquad \uCM \mydef \Exp \left[ u|f \right] = \int u \; p_{post}(u|f) \, \rmd u, \label{eq:PointEst}
\end{equation}
which both yield a single point estimate of $u$. Details on their properties and relationship can be found in \cite{BuLu14}. More sophisticated estimators such as \termabb{conditional covariance}{CCov}, \termabb{conditional variance}{CVar} or \termabb{standard deviation}{CStd} estimates try to extract higher order statistics of $u$, or try to quantify the uncertainties of $u$, for instance through \emph{credible region/interval} and \emph{extreme value probability} estimators. \\
\emph{Bayesian computation} refers to the practical task of computing the above estimators. For most inversion scenarios and prior models, this involves solving high-dimensional optimization or integration tasks (cf., \eref{eq:PointEst}), or even a mix of both. In this article, we are examining techniques that integrate $\post$ by \emph{Monte Carlo integration}:
\begin{equation}
\int g(u) \post \, \rmd u \approx \frac{1}{K} \sum_i^K g(u^i), \label{eq:Quad}
\end{equation}
where $u^i$ are samples of $\post$ generated by a \emph{sampling algorithm}/\emph{sampler}. Due to the lack of efficient \emph{direct samplers} that generate i.i.d. samples, \termabb{Markov chain Monte Carlo}{MCMC} samplers need to be employed in most situations. MCMC for high dimensional Bayesian inversion is a very active field of research, see, e.g., \cite{HaLaMiSa06,CuFoSu11,Ba12,CoRoStWh13,La14,PaMaCh14,SoLuAr14,SoSo14,BaSoHaLa14,AgBaPaSt14,PeScChPeToHeMc15,Pe15,CuLaMa16} for some examples of recent developments.\\
In \cite{Lu12}, an efficient MCMC sampler for Gibbs priors with $\ell_1$-norm-type energies ($\ell_1$-\emph{priors}, cf. Figure \ref{subfig:L1LiPrPo}) was presented:
\begin{equation}
\prior \propto \exp \left(  - \lambda \norm{D^T u}_1 \right). \label{eq:L1Prior}
\end{equation}
Such energies are commonly used to impose \emph{sparsity constraints} on the solution of high dimensional image reconstruction problems, a direction of research closely related to the notion of \emph{compressed sensing} \cite{CaRoTa06,Do06,FoRa13}. A detailed discussion of sparsity as a-priori information in Bayesian inversion can be found in \cite{Lu14}. The sampler developed in \cite{Lu12} belongs to the class of \emph{single component (SC) Gibbs samplers}, which sample $\post$ by subsequently sampling along the conditional, single component densities $p_{post}(u_j | u_{-j},f)$:
\begin{algo}{\textbf{(SC-Gibbs Sampling)}} \label{algo:Gibbs}
Define an initial state $u^0$, a burn-in size $K_0$ and sample size $K$. For $i$ $=$ $1$,$\ldots$,$K_0+K$ do:
 \begin{itemize}[leftmargin=0.15\textwidth  , itemsep=3pt,topsep=5pt]
  \item[A1.1] Choose a component $j$ (deterministic or random). 
  \item[A1.2] Draw  $y \sim  p_{post}( \, \cdot \, | u_{-j},f)$
  \item[A1.3] Set $u^{i+1}_{j} = y$, and $u^{i+1}_{-j} = u^i_{-j}$.
 \end{itemize}
 Discard $\{u^i\}_{i=0}^{K_0}$ and use  $\{u^i\}_{i=K_0+1}^{K_0+K}$ as a sample of $\post$.
\end{algo}
We have used $x_{-j} \mydef (x_1,\ldots,x_{j-1},x_{j+1},\ldots,x_n)^T$ to denote a vector with all components of $x$ except for $x_j$. An illustration of Gibbs sampling is given in Figure \ref{subfig:Gibbs}. In \cite{Lu12}, this SC-Gibbs sampler was compared to the popular \termabb{Metropolis-Hastings }{MH} sampler: For the computational scenarios considered and the evaluation performed, it was demonstrated that in contrast to the MH, SC-Gibbs sampling gets more efficient when the level of sparsity or the dimension of the unknowns is increased. Thereby, it became possible to carry out sample-based inference with $\ell_1$ priors in challenging inverse problems scenarios with $n > 10^6$:
\begin{itemize}
\item The theoretical predictions about the infinite dimensional limits of TV priors posed in \cite{LaSi04,Co11} could be verified numerically (see \cite{Lu14}). 
\item \termabb{Computed tomography}{CT} inversion with Besov space priors (cf. \cite{KoLaNiSi12,HmKaKoLaNiSi13}) was examined for simulated and experimental data (see \cite{BuLu14,Lu14}). 
\item The numerical results stimulated the development of new theoretical ideas about the relationship of MAP and CM estimates (see \cite{BuLu14,HeBu15}).
\end{itemize}

\subsection{Previous Limitations} \label{subsec:PreviousLimitations}

As the sampler developed in \cite{Lu12} relies on a direct sampling of the SC densities, namely the \termabb{inverse cumulative distribution method}{iCDF}, we will call it the \emph{direct $\ell_1$ sampler} from now on. While the direct $\ell_1$ sampler works well in the applications described above, it suffers from several limitations. To understand them, we recall that an efficient SC-Gibbs sampler needs to 
\begin{itemize}[leftmargin=0.15\textwidth  , itemsep=3pt,topsep=5pt]
   \item[(SC1)] compute the conditional, SC densities in an explicit, parameterized form in a fast way.  
   \item[(SC2)] employ a fast, robust and exact sampling scheme for the parameterized form of the SC densities. \label{item:SC2} 
\end{itemize}
In order to best fulfill (SC1) and (SC2), the direct $\ell_1$ sampler was designed for a very particular setting: Firstly, in addition to relying on a linear forward map \eref{eq:FwdEq} and a Gaussian noise model \eref{eq:Likelihood}, it assumes that the operator $D \in \R^{n \times h}$ in \eref{eq:L1Prior} can be diagonalized (\emph{synthesis} priors): There is a basis matrix $V$ such that $D^T V$ is a diagonal matrix $W \in \R^{h \times n}$. The direct $\ell_1$ sampler then samples over the coefficients of $u = V \xi$ in this basis:
\begin{equation}
p_{post}(\xi|f) \propto \exp \left( -\case{1}{2} \sqnorm{f - AV \, \xi}  - \lambda \, \norm{ W \xi}_1 \right)
\end{equation}
This excludes the use of \emph{frames} or \emph{dictionaries} for $D$. Secondly, it only works for the $\ell_1$ norm as a prior energy: A straight-forward extension of iCDF to examine more general $\ell_p^q$-prior of the form $\prior \propto \exp \left( - \lambda \norm{ D^T u }_p^q \right)$ is not possible. This excludes the interesting cases of $q = p, p < 1$, which leads to a non-convex energy but also $p = 1$, $q > 1$, which was examined in \cite{Co11}. Finally, a lot of interesting priors such as the popular isotropic TV prior in 2D/3D or related, \emph{bloc/structured sparsity} priors have a more involved structure than \eref{eq:L1Prior} and cannot be treated with iCDF in an efficient and robust way as well. In all the above cases, including additional \emph{hard constraints}, $u \in \mathcal{C}$, where $\mathcal{C}$ is the \emph{feasible set} of solutions is often advantageous:
\begin{equation}
\prior \propto \tilde{p}_{prior}(u) \cdot \Indicator_\mathcal{C} (u) = \cases{
\tilde{p}_{prior}(u) & if $u \in \mathcal{C}$ \\
0	& otherwise} \label{eq:HardCon}
\end{equation}
While such constraints have proven to be very useful as a-priori knowledge \cite{Vo02,BaFo12}, their inclusion into the direct $\ell_1$ sampler in a numerically stable way is cumbersome. 

\subsection{Contributions and Structure} \label{subsec:ConStr}

For most of the limitations discussed above, the main problem is not to fulfill (SC1), but to fulfill (SC2) by using a direct sampler such as iCDF for the parameterized SC densities in step A1.2. In this article, we sample from them by using a generalization of \emph{slice sampling} that utilizes their specific structure instead and demonstrate the effectiveness of this replacement in different computed examples. This allows us to perform sample-based Bayesian inference in high-dimensional scenarios with the priors described above for the first time. \\
The paper is structured as follows: In Section~\ref{sec:Methods}, we first derive the SC densities for the priors discussed above. Then, we introduce the basic and generalized slice sampler and discuss how to integrate it into the SC Gibbs sampler for Bayesian inversion. Section~\ref{sec:CompExa} contains computed examples and Section~\ref{sec:Discussion} closes with a discussion. Several technical details are covered in Section~\ref{sec:Impl}.


\section{Sampling Methods} \label{sec:Methods}

For general and comprehensive introductions to MCMC sampling methods, we refer to \cite{RoCa05,Li08}.

\subsection{SC Posteroir Densities} \label{subsec:SCdensities}

In this section, we briefly derive the SC posterior densities for the examined prior models in a simple, parameterized way, cf. (SC1). We first discuss the case where a basis $\{v_1,\ldots,v_n\}$ helps to represent $u = \sum \xi_l v_l \mydefswitched V \xi$ such that $p_{post}(\xi_j | \xi_{-j}, f)$ can be described using as few parameters as possible. Once such a basis is found, the part of $p_{post}(\xi_j | \xi_{-j}, f)$ coming from the likelihood is easy to derive: We define $\Psi \mydef AV$ and $\varphi(j) \mydef f - \Psi_{-j} \xi_{-j}$. Then, we find that
\begin{eqnarray}
\fl \quad  \case{1}{2} \sqnorm{f - Au}   =  \case{1}{2} \sqnorm{f- A V \xi} = \case{1}{2} \sqnorm{f- \Psi \xi}  = 
\case{1}{2} \sqnorm{ f - (\Psi_{-j}  \xi_{-j} + \Psi_j \xi_j)} \nonumber \\ 
\qquad  = \case{1}{2} \sqnorm{ \varphi(j) - \Psi_j \xi_j} 
\proptoab{\xi_j}  \case{1}{2} \sqnorm{\Psi_j} \xi_j^2 + \Psi_j^T \varphi(i) \xi_j \mydefswitched  a x^2 - bx,
\end{eqnarray}
where we introduced $x \mydef \xi_j$, $a \mydef \case{1}{2} \sqnorm{\Psi_j}$, and $b \mydef \Psi_j^T \varphi(j) = \Psi_j^T f - (\Psi_j^T \Psi_{-j}) \xi_{-j}$ to ease the notation for the following sections. Note that while $a$ and $ \Psi_j^T f$ can be precomputed, $(\Psi_j^T \Psi_{-j}) \xi_{-j}$ relies on the current state of the $\xi$-chain and has to be computed in every step of the sampler. Especially for complicated forward operators in high-dimensional scenarios, this operation is the computational bottleneck of SC Gibbs samplers. Therefore, a careful, scenario-dependent implementation is important to obtain a fast sampler.  \\
Now we proceed to determine $V$ and the part of $p_{post}(\xi_j | \xi_{-j}, f)$ coming from the prior. The energies of the $\ell_p^q$ priors can be written as
\begin{equation}
\J(u) = \left( \sum_k^h \abs{D_k^T u }^p  \right)^{\frac{q}{p}} =  \left( \sum_k^h   \Big| \sum_l (D_k^T v_l) \xi_l \Big|^p \right)^{\frac{q}{p}}. \label{eq:J1} 
\end{equation}
To obtain simple conditional densities for all $\xi_j$, we thus have to choose $V$ such that 
\begin{equation}
\max_l \abs{ D^T v_l}_0, \qquad \text{where} \qquad \abs{ u }_0 \mydef \text{card}\left( \{ k | u_k \neq 0 \} \right),\\
\end{equation}
is as small as possible. We first consider the special but important case of $D^T \in \R^{h \times n}$ having full rank and $h \leqslant n$. This includes the case where the columns $D$ are elements of a basis, and thereby, the class of \emph{Besov priors}, see \cite{LaSaSi09,KoLaNiSi12,HmKaKoLaNiSi13,DaHaSt12,BuLu14} and the TV prior in 1D with Neumann boundary conditions, which we will use in the computational studies. Due to the full rank, we can choose $v_1,\ldots,v_h$ such that $D^T v_l = e_l$ for $l = 1,\ldots,h$, and $v_{h+1},\ldots,v_n$ such that $D^T v_l = 0$ for $l = h+1,\ldots,n$ (for $D$ being a basis, we have $V = D$). With this transformation, \eref{eq:J1}  simplifies to
\begin{equation}
\J(\xi) \propto  \left( \sum_l^h \abs{\xi_l}^p \right)^{\frac{q}{p}}   = \left(  \abs{\xi_j}^p + \sum_{l \neq j}^h  \abs{\xi_l}^p \right)^{\frac{q}{p}}.
\end{equation}
Defining $x \mydef \xi_j$ as above, we can write the conditional SC posterior density as
\begin{equation}
\fl \quad p(x) \propto \exp \left(- a x^2 + b x - c \left(|x|^p + d \right)^{q/p} \right), \quad c \mydef \lambda\, \Indicator_{\{j \leqslant h \}}, \quad d \mydef \sum_{l \neq j}^h  \abs{\xi_l}^p, \label{eq:SCLpq}
\end{equation}
 which simplifies to 
 \begin{equation}
p(x) \propto \exp \left(- a x^2 + b x - c |x|^p\right), \quad c \mydef \lambda \, \Indicator_{\{j \leqslant h \}}, \label{eq:SCLp}
\end{equation}
for $\ell^p_p$ priors. In the case where $D$ cannot be diagonalized, an explicit form is given by 
 \begin{eqnarray}
\fl \qquad \qquad p(x) \propto \exp \left(- a x^2 + b x - c \left( \sum_{k \in \text{supp}{(D^T v_j)}} \abs{d_k x - e_k }^p \right)^{\frac{q}{p}}\right), \label{eq:SCgenLpq}\\
\fl \qquad \text{where}  \qquad  c \mydef \lambda \, \Indicator_{\{j \leqslant h \}}, \qquad d_k \mydef \left( D^T v_{j} \right)_k, \qquad e_k \mydef \left( D^T V_{-j} \xi_{-j}\right)_k.
\end{eqnarray}
Various generalizations of the standard $\ell_p^q$ priors with $\norm{D^T u}_p^q$-type energies first compute the $\ell_2$-norm of a \emph{local feature} of $u$, e.g., of its gradient, and then measure the \emph{global} $\ell_p^q$ energy of these local $\ell_2$ norms. In this article, we will only discuss one prominent example thereof, which is given by the isotropic TV prior in 2D: If we assume that $u$ represents an $N \times N$ discrete image, we can index the components of $u$ as $u_{(k,l)}$ with $k=1,\ldots,N$, $l=1,\ldots,N$, $n = N^2$. We can then use forward differences in both spatial directions to define
\begin{equation}
\fl \qquad \J_{iTV}(u) = \sum_{(k,l)}^n \sqrt{(u_{(k+1,l)}-u_{(k,l)})^2 + (u_{(k,l+1)} - u_{(k,l)})^2}, \label{eq:iTVprior}
\end{equation}
with appropriate additional boundary conditions. The local nature of the $\J_{iTV}(u)$ allows to derive a simple parameterization of the SC densities in the pixel basis $V = I_n$. Every $\xi_j = u_{(k,l)}$ only appears in three terms of the energy:
\begin{eqnarray}
\fl \qquad \J_{iTV}\left( u_{(k,l)} \, \big| \, u_{-(k,l)}\right) \proptoab{(k,l)}  &\phantom{,...} \sqrt{(u_{(k+1,l)}-u_{(k,l)})^2 + (u_{(k,l+1)} - u_{(k,l)})^2 } \nonumber \\
&+ \sqrt{(u_{(k,l)}-u_{(k-1,l)})^2 + (u_{(k-1,l+1)} - u_{(k-1,l)})^2 } \nonumber \\
&+ \sqrt{(u_{(k+1,l-1)}-u_{(k,l-1)})^2 + (u_{(k,l)} - u_{(k,l-1)})^2 }
\end{eqnarray}
Therefore, we can write the conditional SC posterior as 
\begin{equation}
\quad p(x) \propto \exp \left( -a x^2 + bx - c \sum_{k=1}^3 \sqrt{ d_k (x - e_k)^2  + g_k }\right), \label{eq:SCTV}
\end{equation}
with appropriately computed parameters $d_k \in \{0,1,2\}$, $e_k$, $g_k \geqslant 0$. \\
The difficulty of incorporating additional hard constraints \eref{eq:HardCon} depends on the shape of the feasible set $\mathcal{C}$ and the transformation $V$ applied. In the following, we assume that they lead to a feasible (semi-)finite interval $[lb,ub]$ to which the continuous densities computed above can be restricted to. In the case of $\mathcal{C}$ being convex, such an interval always exists and there are computationally efficient ways to compute it.

\subsection{MCMC-within-Gibbs Sampling} \label{subsec:MCMCwithinGibbs}

The direct $\ell_1$ sampler is sampling \eref{eq:SCLp} with $p = 1$ by the iCDF method using an explicit form of the inverse CDF. For $p,q \neq 1$ or \eref{eq:SCTV} this is not possible and one would need to integrate the CDF numerically to use the iCDF method as a SC density sampler. However, already for $p = 1$, a major technical difficulty was to develop a numerical implementation that worked for all possible combinations $(a,b,c)$: The first implementations broke down when the dimension $n$ of the problem was increased and the ill-posedness became more severe. The reason was that combinations of $(a,b,c)$ corresponding to extremely degenerate SC densities appeared more frequently for $n \rightarrow \infty$ and in general, the variability of SC densities grows. This trend will be an even more severe problem when one cannot find an explicit form of the inverse CDF and needs so resort to numerical integration. But also replacing the iCDF method by a univariate MCMC sampler (\emph{MCMC-within-Gibbs} sampling) becomes challenging: The most commonly used \emph{Metropolis-within-Gibbs} sampler, which utilizes an easy-to-implement MH sampler with a univariate Gaussian proposal $\mathcal{N}(x,\kappa)$ (where $x$ is the current state) for the SC sampling step A1.2 will not work properly in such a situation: For an MH sampler to be efficient, finding a value of $\kappa$ leading to an optimal \emph{acceptance rate} is essential. However, the large variations in-between SC densities renders an automatic tuning of a single $\kappa$ impossible. The alternative would be to tune and use a different $\kappa_j$ for every component $j$, but the tuning procedure would require $n$ times more samples than tuning one $\kappa$ for all components. Thereby, the resulting algorithm would be more like an \emph{adaptive} SC-MH sampler than a Gibbs sampler \cite{HaSaTa05,LaRoRo13}.

\subsection{Slice Sampling} \label{subsec:SliceSampler}

\begin{figure}[tb]
\centering
\subfloat[][]{\includegraphics[width = 0.3215\textwidth]{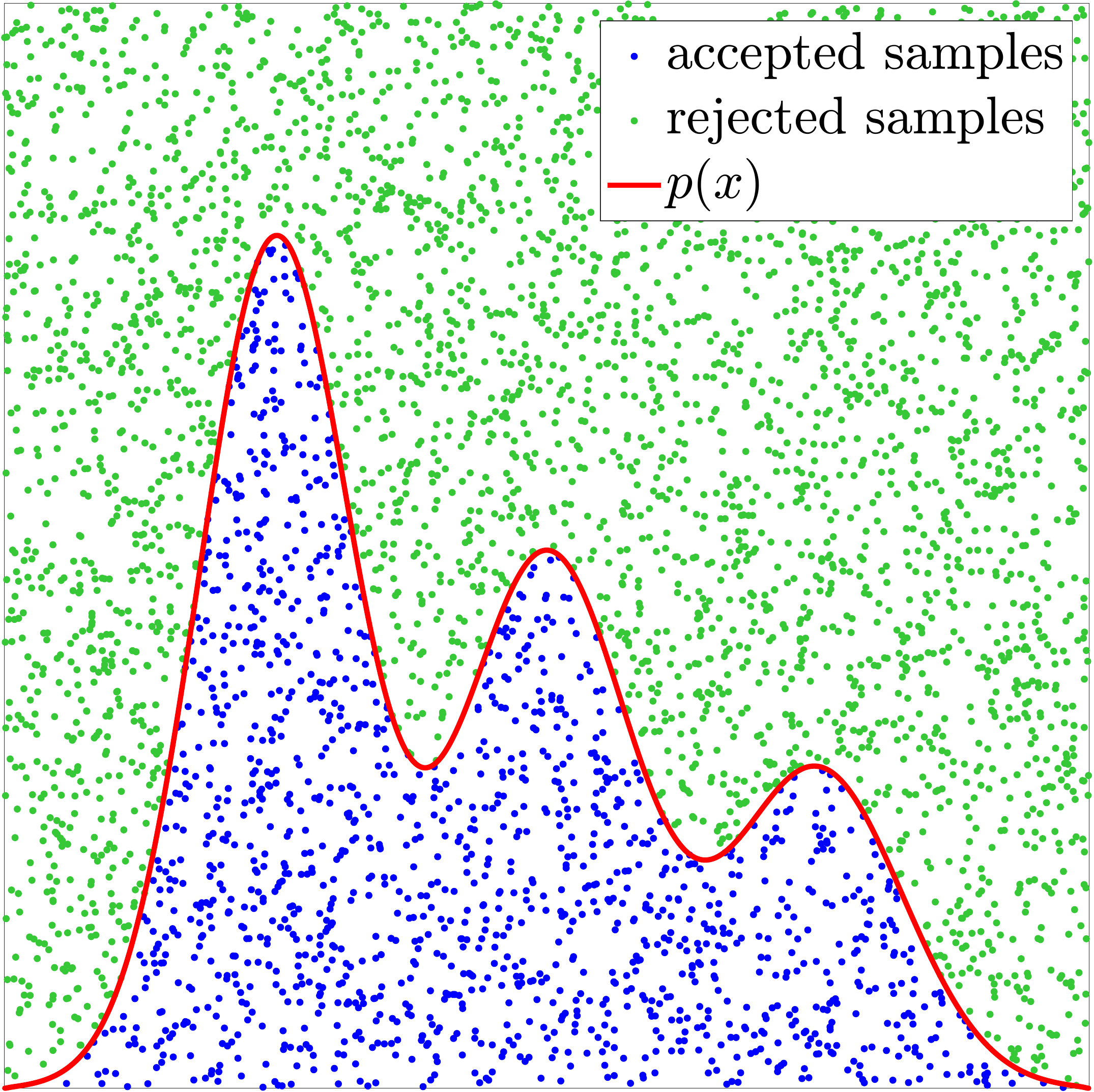} \label{subfig:RejectSamplingA}}
\hskip 0.005\textwidth
\subfloat[][]{\includegraphics[width = 0.32\textwidth]{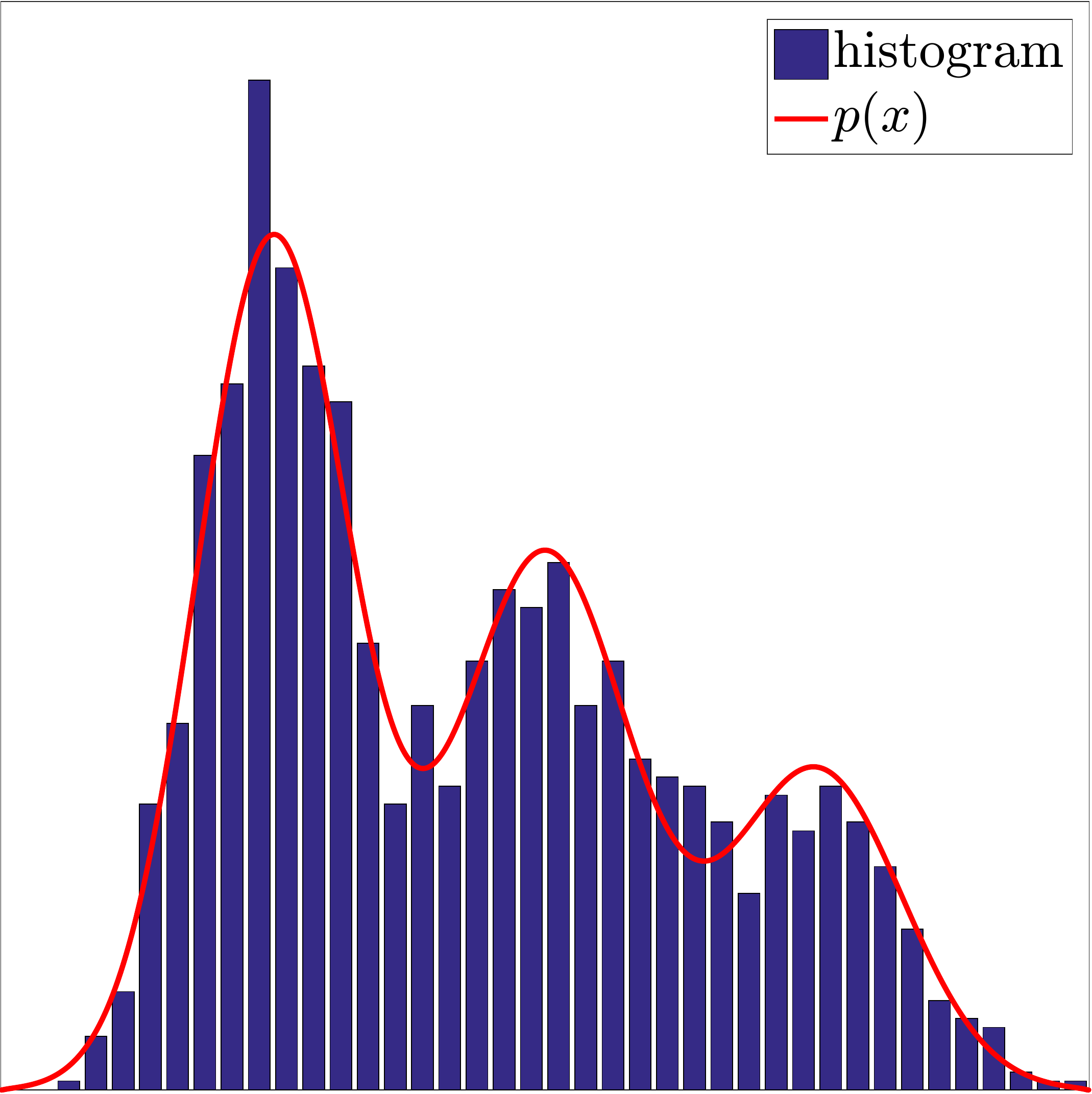} \label{subfig:RejectSamplingB}}
\hskip 0.005\textwidth
\subfloat[][]{\includegraphics[width = 0.32\textwidth]{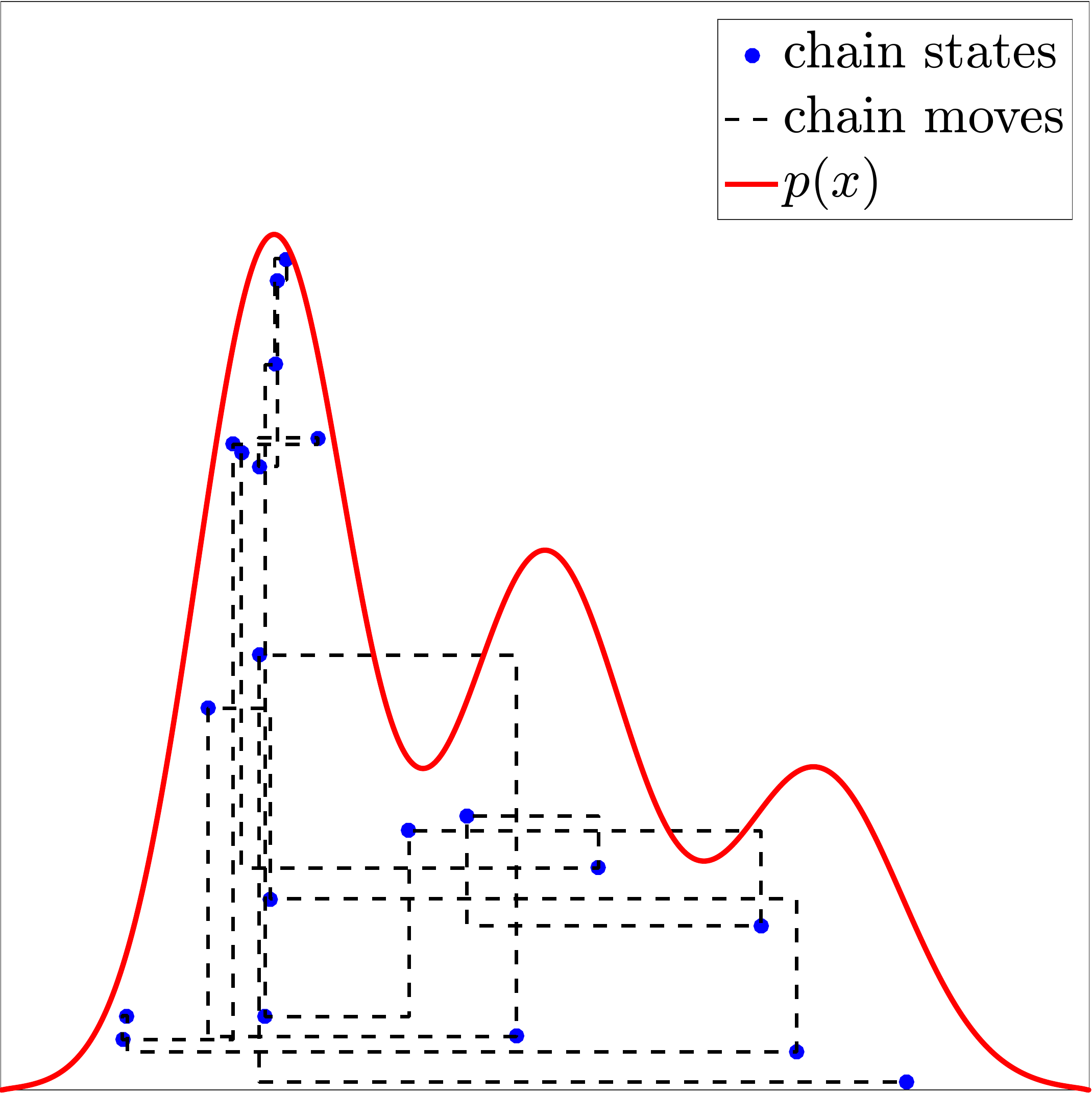} \label{subfig:SliceSampler}}
\caption{An illustration of accept-reject methods and slice sampling: \protect\subref{subfig:RejectSamplingA} To sample from the density $p(x)$ (red line), uniform samples $(x^i,y^i)$ (blue and green dots) are generated in a region enclosing its graph. All samples fulfilling $y^i \leqslant p(x^i)$ (blue dots) are accepted. \protect\subref{subfig:RejectSamplingB} Histogram computed from the $x$ values of all accepted samples. \protect\subref{subfig:SliceSampler} Slice sampling. The $x$ coordinates of the blue dots are samples of $p(x)$, while the dashed black line illustrates the path of the sampler on $\mathcal{G}_p$. \label{fig:RejectSliceSampling}}
\end{figure}

\emph{Slice sampling} transfers the automatic adaptation of Gibbs sampling to univariate densities. While the basic version to sample arbitrary densities in a "back-box" fashion was proposed in \cite{Ne03}, we follow the presentation given in \cite{RoCa05}, which leads to a general version in which we can utilize several properties of our specific posterior densities to derive an efficient SC sampler. The starting point for slice sampling is the \emph{Fundamental Theorem of Simulation}, which states that sampling from a distribution $p(x)$ is equivalent to sampling uniformly from the area under the graph of $p(x)$: $\mathcal{G}_p \mydef \{(x,z)| 0 \leqslant z \leqslant p(x) \}$. This simple observation is the basis of \emph{accept-reject samplers}, a widely used class of samplers which draw uniform samples $(x,z)$ from a region enclosing $\mathcal{G}_p$ and only accept the sample if it fulfills $z \leqslant p(x)$. Figures \ref{subfig:RejectSamplingA} and \ref{subfig:RejectSamplingB} illustrate this principle. Slice sampling utilizes this principle in another way: It samples the auxiliary, bivariate density $\tilde{p}(x,z) \propto  \Indicator_{\mathcal{G}_p}(x,z)$ by a Gibbs sampler and only keeps the $x$ samples, cf. Figure \ref{subfig:SliceSampler}:
\begin{algo}{\textbf{(Basic Slice Sampling)}} \label{algo:BasicSlice}
For a univariate density $p(x)$, define an initial state $x^0$, a burn-in size $K_0$ and a sample size $K$. For $i$ $=$ $1$,$\ldots$,$K_0+K$ do:
 \begin{itemize}[leftmargin=0.15\textwidth  , itemsep=3pt,topsep=5pt]
  \item[A2.1] Draw $y$ uniform from $\left[0,p(x^i)\right]$ (vertical move).
  \item[A2.2] Draw $x^{i+1}$ uniform from $S^y \mydef \{z \;|\; p(z) \geqslant y \}$ (horizontal move).
 \end{itemize}
Discard $\{x^i\}_{i=0}^{K_0}$ and use  $\{x^i\}_{i=K_0+1}^{K_0+K}$ as a sample of $p(x)$.
\end{algo}
The difficulty of this basic slice sampling scheme as developed in \cite{Ne03} is determining $S^y$ in Step A2.2. For the SC densities we want to sample from, determining $S^y$ explicitly is not always feasible, and robust numerical approaches to compute it are difficult to design. For instance, using non-convex prior energies such as in $\ell_p^q$ priors with $q = p, p < 1$ leads to multi-modal SC densities and $S_y$ may not be a single interval. Therefore, we will use a generalization of Algorithm \ref{algo:BasicSlice}: Slice sampling is a variant of \emph{auxiliary variables algorithms} that introduce an additional variable $y$ with a suitable density $p(y|x)$. Then, samples $(x^i,y^i)$ from $p(x,y) = p(x) p(y|x)$ are obtained by a Gibbs sampler, which relies on $p(y|x)$ and $p(x|y)$, and only the $x^i$ are kept. For the basic slice sampler, $p(y|x)$ is chosen as 
\begin{equation}
p(y|x) = \frac{1}{p(x)} \Indicator_{\{[0,p(x)]\}}(y),
\end{equation}
i.e., as a uniform distribution on $[0,p(x^i)]$. We then have
\begin{eqnarray}
p(x,y) &= p(x) \frac{1}{p(x)} \Indicator_{[0,p(x)]}(y)\\
p(x|y) &\propto \Indicator_{\{[0,p(x)]\}}(y) = \Indicator_{\{x \;|\; p(x) \geqslant y\}}(x)
\end{eqnarray}
If $p(x)$ factorizes to $p(x) \propto p_1(x)  p_2(x)$ we can define
\begin{equation}
p(y|x) = \frac{1}{p_2(x)} \Indicator_{\{[0,p_2(x)]\}}(y),
\end{equation}
which leads to
\begin{eqnarray}
\fl \quad  &p(x,y) = p(x) p(y|x) = p(x) \frac{1}{p_2(x)} \Indicator_{\{[0,p_2(x)]\}}(y) = p_1(x) \Indicator_{\{[0,p_2(x)]\}}(y),\\
\fl \quad  &p(x|y) = p_1(x)  \Indicator_{S_2^y}(x), \quad \text{with} \quad S_2^y \mydef \left\lbrace z \;|\; p_2(z) \geqslant y \right\rbrace \label{eq:S2Set}
\end{eqnarray}
The corresponding sampler takes the form:
\begin{algo}{\textbf{(Slice Sampling)}\\} \label{algo:Slice}
For a univariate density $p(x) \propto p_1(x) p_2(x)$, define an initial state $x^0$, a burn-in size $K_0$ and sample size $K$. For $i$ $=$ $1$,$\ldots$,$K_0+K$ do:
 \begin{itemize}[leftmargin=0.15\textwidth, itemsep=3pt,topsep=5pt]
  \item[A3.1.] Draw $y$ uniform from $\left[0,p_2(x^i)\right]$ (vertical move).
  \item[A3.2.] Draw $x^{i+1}$ from $p_1(x) \Indicator_{S^y_2}(x)$ (weighted horizontal move).
 \end{itemize}
 Discard $\{x^i\}_{i=0}^{K_0}$ and use  $\{x^i\}_{i=K_0+1}^{K_0+K}$ as a sample of $p(x)$.
\end{algo}
For all the methods presented in this section, $p(x)$ does not need to be normalized. Also note that for simplicity, we refer to Algorithm \ref{algo:Slice} as the "slice sampler", hopefully without causing confusion with the one presented in \cite{Ne03}, which was included as the "basic slice sampler" (Algorithm \ref{algo:BasicSlice}) here for completeness of the presentation.


\subsection{Slice-Within-Gibbs Sampling for Bayesian Inversion} \label{subsec:SWG}

The implicit variable split introduced in Algorithm \ref{algo:Slice} is appealing if $S_2^y = \left\lbrace z \;|\; p_2(z) \geqslant y \right\rbrace$ is a single interval and easy to determine and $p_1(x)$ constrained to an interval is easy to sample from. For the SC posterior densities we consider here, this holds if we split into likelihood plus hard constraints, i.e., $p_1(x) = \exp(-ax^2 + bx) \Indicator_{[lb,ub]}(x)$, and prior parts $p_2(x)$. As the prior terms are unimodal and some even symmetric to zero, $S_2^y$ is a single interval and can be determined easily: For  \eref{eq:SCLpq}, we have $p_2(x) \propto  \exp\left(-c\left(|x|^p+d\right)^{q/p}\right)$ and
\begin{equation}
\exp\left(-c\left(|x|^p+d\right)^{q/p}\right) \geqslant y \Longleftrightarrow |x| \leqslant \left(\left( -\frac{\log(y)}{c}\right)^{p/q} - d \right)^{1/p}.
\end{equation}
For the TV prior, \eref{eq:SCTV}, we need to compute $S_2^y$ numerically. However, as the energy of $p_2(x)$ is convex,  $S_2^y$ is a single interval given by the solutions to $p_2(x) = y$. As the energy of $p_2(x)$ is also piecewise smooth and can be bounded from below, we can easily find starting points for fast, derivative-based root-finding-algorithms. The details are given in \ref{subsec:ImplTVSS}. A generalization to other convex, piecewise-smooth energies, such as \eref{eq:SCgenLpq} with suitable $p, q$, is straight-forward ($p = q = 1$ is a special case as $p_2(x) = y$ can be solved explicitly by a simple scheme). However, if $D^T V$ is dense the number terms in the prior energy is large and this step will become the computational bottleneck of the whole solver. Fortunately, many relevant operators $D^T$ such as finite difference operators or dictionaries composed of local patches are sparse in the original basis, $V = I_n$. \\
The likelihood part $p_1(x)$ is a Gaussian with $\mu_{SS} = b/(2a)$ and $\sigma^2_{SS} = 1/(2a)$, truncated to the interval $I = S_2^y \cap [lb,ub]$. For sampling truncated Gaussians, various direct samplers were developed. Our implementation relies on a modified, more robust, version of \cite{Ch12}. Note that if the sampler is initialized in a feasible point $u^0 \in \mathcal{C}$, the probability of $I$ being empty or a single point is zero in theory. In practice, finite precision can lead to $I = \{\tilde{x}\}$, in which case one has to set $x^{i+1} = \tilde{x}$. \\
Using the slice sampler presented above to sample from $p_{post}( \, \cdot \, | \xi_{-j},f)$ in step A1.2 will be called \emph{slice-within-Gibbs} sampling. In principle, it will generate a full Markov chain 
\begin{equation}
\left\lbrace \xi_j^{i_s} \right\rbrace_{i_s=K_{s,0}+1}^{K_{s,0}+K_s} ~ \sim p_{post}( \, \cdot \, | \xi_{-j},f),
\end{equation}
where we subscripted all variables belonging to the inner slice sampler with $s$. Practically, we only need one sample from $p_{post}(  \, \cdot \, | \xi_{-j},f)$. We will always initialize the slice sampler with the current value $\xi_i$ of the component we want to update. Then, we only have to determine the length of the burn-in phase $K_{s,0}$ and choose the first sample of the real run as a sample of $p_{post}( \, \cdot \, | \xi_{-j},f)$, i.e., $K_s = 1$.\\
The correctness and convergence of the slice-within-Gibbs sampler can be established by combining the properties of the slice sampler (Algorithm \ref{algo:Slice}) and the general Gibbs sampler (Algorithm \ref{algo:Gibbs}), which are discussed in \cite{RoCa05}.


\section{Computed Examples} \label{sec:CompExa}

\subsection{Computational Scenarios}

\begin{figure}[bt]
   \centering
\subfloat[][Unknown function \label{subfig:BoxcarReal}]{\includegraphics[height = 0.25\textwidth]{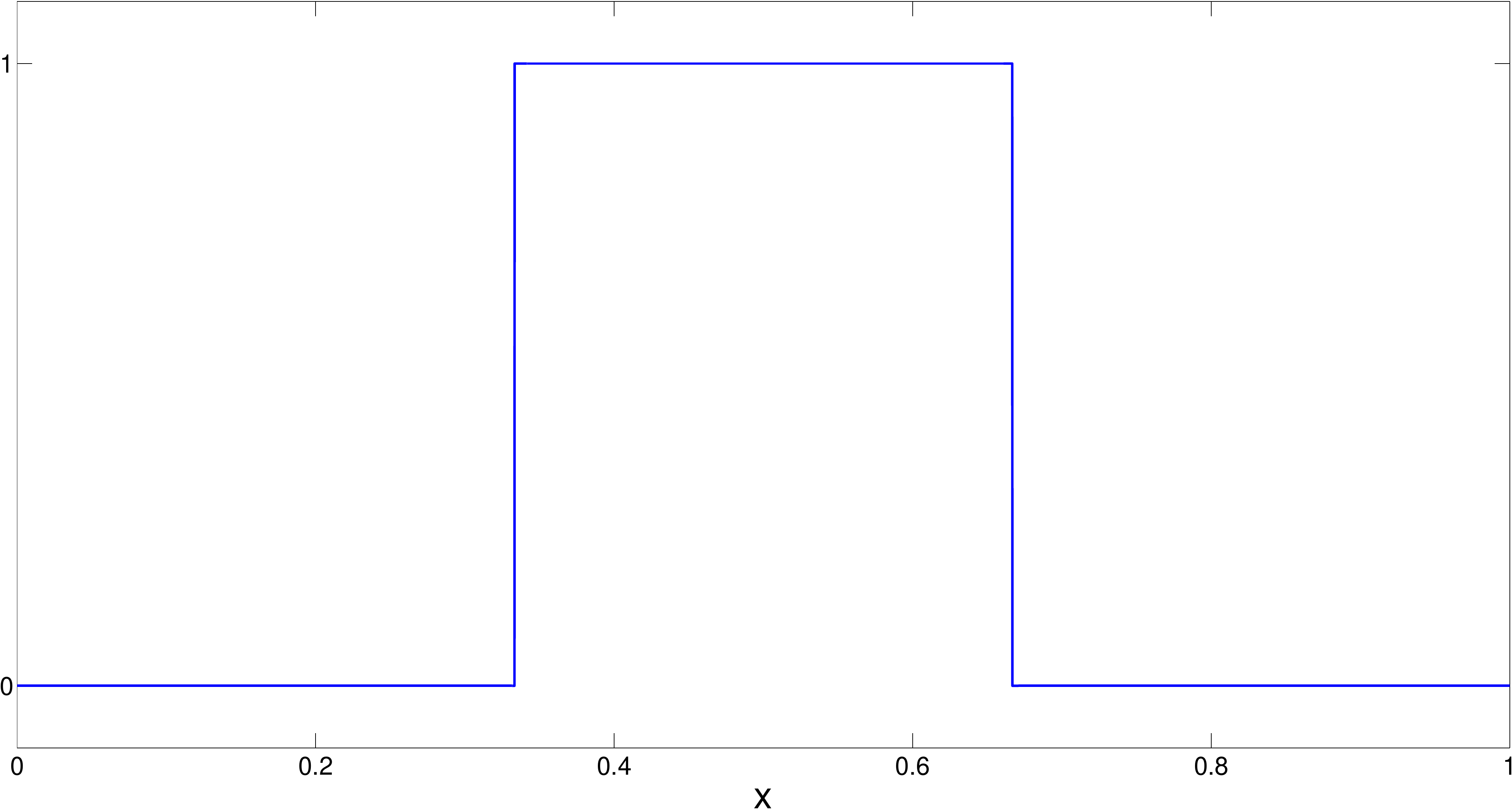}}
\hspace{0.01\textwidth}
\subfloat[][Measured data $f$ \label{subfig:BoxcarDataClean}]{\includegraphics[height = 0.25\textwidth]{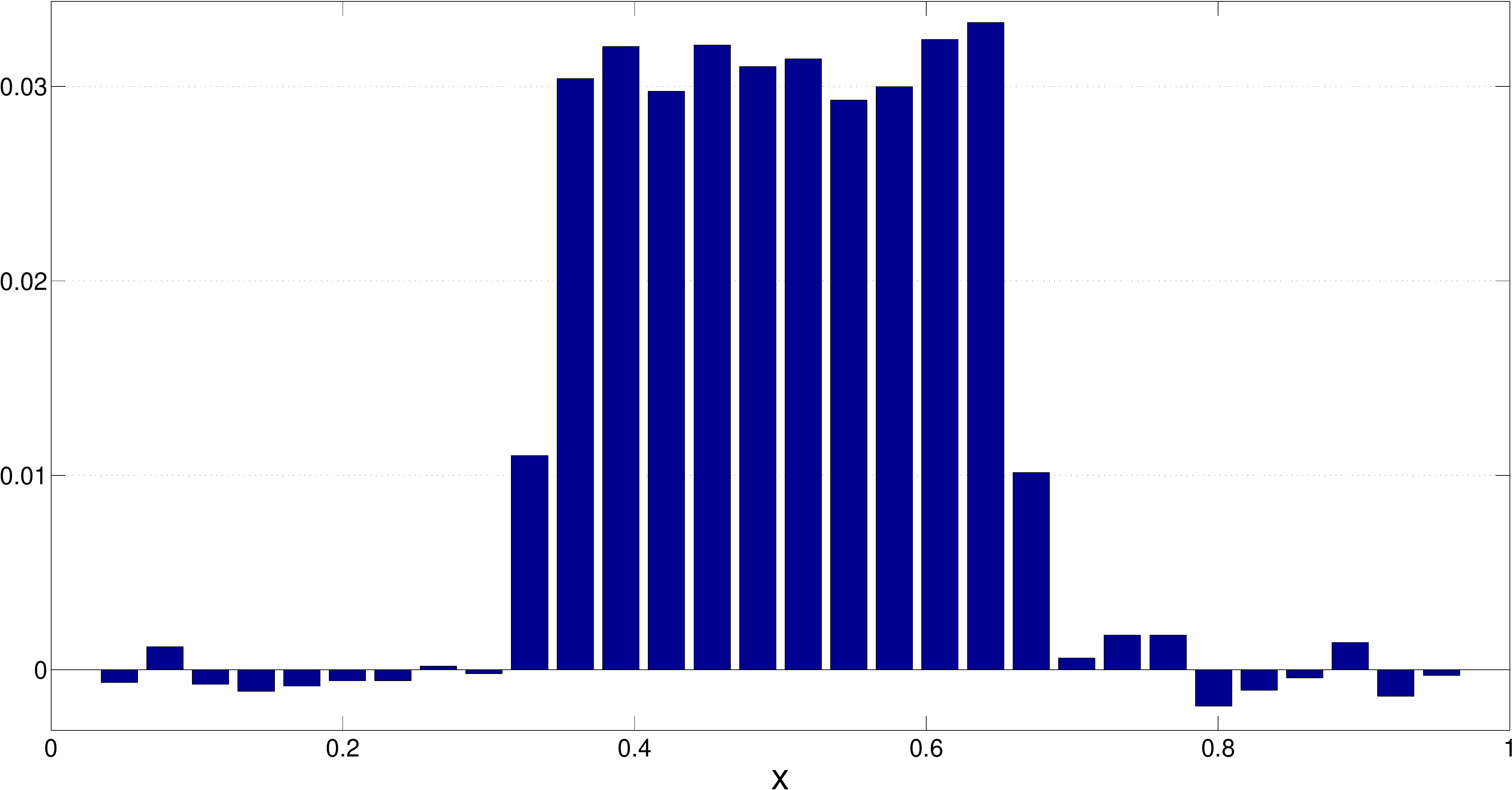}}
\caption{One dimensional image deblurring scenario ("Boxcar").}
   \label{fig:BoxcarSce}
\end{figure}

\subsubsection*{"Boxcar" - Image Deblurring in 1D}
For the initial evaluation studies, we use a simple image deblurring scenario in 1D that was adopted from \cite{LaSi04} and also used in \cite{Lu12}. It is a simplification of  the task to reconstruct a spatially distributed intensity image that is known to consist of piecewise homogeneous parts with sharp edges: The indicator function of $[\case{1}{3},\case{2}{3}]$  is to be recovered from its integrals over $m = 30$ equidistant subintervals of $[0,1]$, corrupted by noise with $\mu = 0, \Sigma = 10^{-3} I_m$ (see Figure \ref{fig:BoxcarSce}). The reconstruction is carried out on the grid $x^n_i = \frac{i}{256}$, $i = 1,\ldots,n$, with $n = 255$ and the forward operator is discretized by the trapezoidal quadrature rule applied to that grid. Further details can be found in Section 3.1.1 of \cite{Lu12}. \\
The prior operator $D^T$ will be given by the forward difference operator with Neumann boundary conditions:
\begin{equation}
\fl \qquad D_i = e_{i+1} - e_i, \qquad \Longrightarrow D_i^T u =  u_{i+1} - u_i, \qquad i=1,\ldots, n-1 \label{eq:IncD}
\end{equation}  
$D^T$ has full rank $h = n-1$ and $V \in R^{n \times n}$ given by 
\begin{equation}
 V_{(i,j)} = \cases{1 & if $i \geqslant j$\\
	    0 & else}
\end{equation}
is a basis matrix $V$ such that $D^T V$ is a diagonal matrix. We will refer to priors based on this operator as \emph{increment priors}. For the $\ell_1$ increment prior, i.e., the conventional TV prior, we can also use the direct $\ell_1$ sampler to sample from the posterior. By this, we can validate the approximation of the direct SC sampling via iCDF by the slice sampler proposed here. We will refer to this setting as the "Boxcar" scenario in the following.

\subsubsection*{"Phantom-CT" - CT Inversion in 2D}

\begin{figure}[tb]
\centering
\subfloat[][\label{subfig:PhanCTReal}]{\includegraphics[height = 0.5\textwidth]{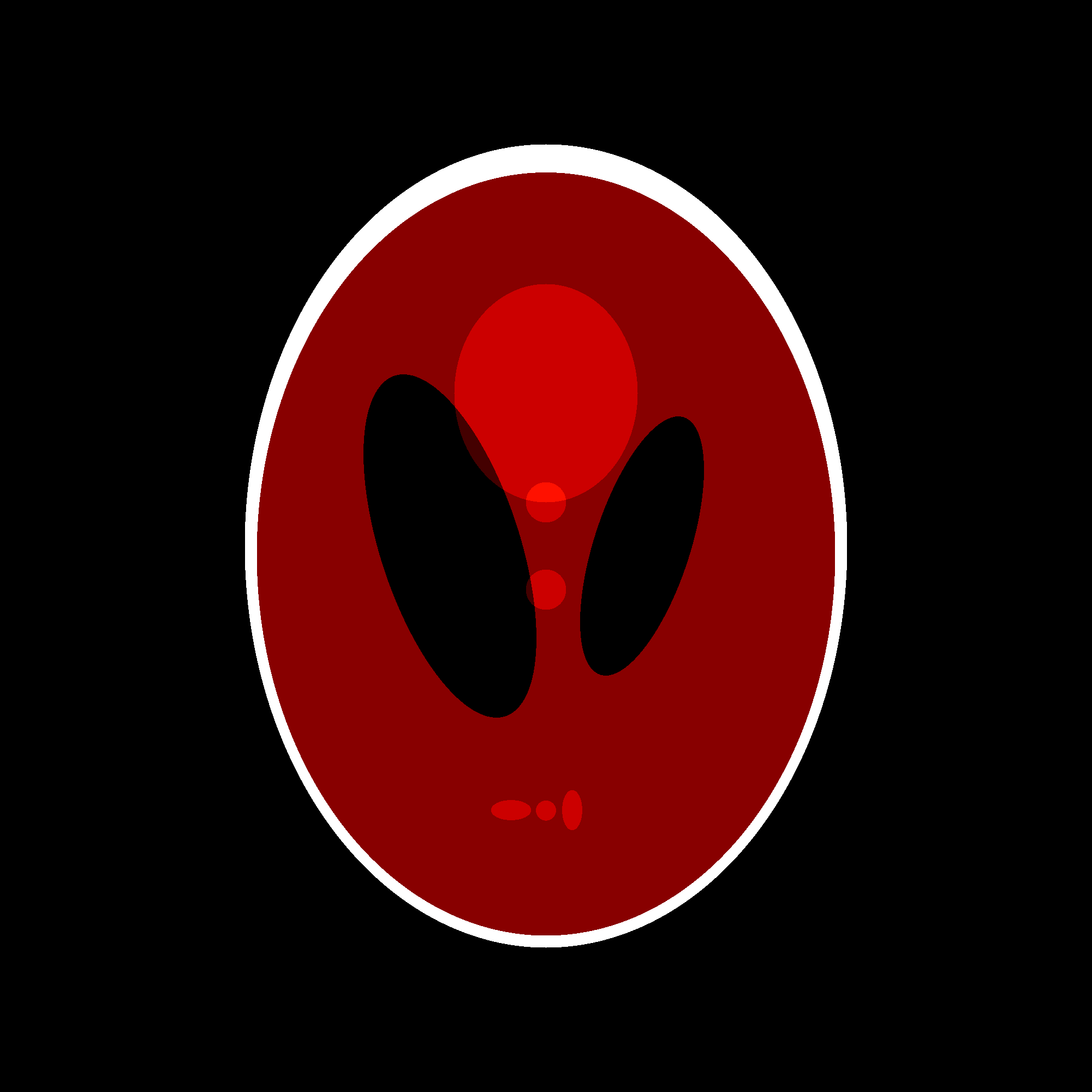}}
\hspace{0.02\textwidth}
\subfloat[][\label{subfig:PhanCTSparseData}]{\includegraphics[height = 0.5\textwidth]{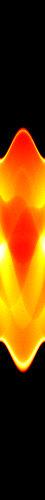}}
\hspace{0.02\textwidth}
\subfloat[][\label{subfig:PhanCTColormap}]{\fbox{\includegraphics[height = 0.5\textwidth]{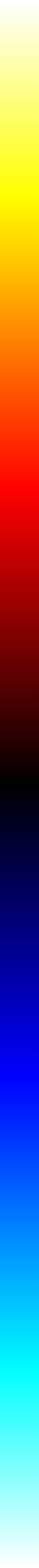}}}
\caption[]{"Phantom-CT" scenario. \subref{subfig:PhanCTReal} Unknown function to recover. 
\subref{subfig:PhanCTSparseData} Clean measurement data ("sinogram") for $m_s = 500$, $m_\theta = 45$. \subref{subfig:PhanCTColormap} Colormap used for all visulizations in this scenario, $0$ corresponds to black. \label{fig:PhanCT}}
\end{figure}

We consider an example of 2D \emph{sparse angle} CT to demonstrate the potential of the proposed sampler for real-world applications. An approximate model of CT is given by the \emph{Radon transform} $\Radon$: For a 2D function $\uInf \in L_2([-1,1]^2)$, it computes integrals along straight lines which are parametrized by the angle $\theta$ of their normal vector and their (signed) distance $s$ to the origin:
\begin{eqnarray}
\Radon[\uInf] (\theta, s) &= \int_{l(\theta,s)} \uInf(x(t),y(t)) \, \rmd l(t) \nonumber \\
&= \int_{-\infty}^{\infty} \uInf( t \sin \theta + s \cos \theta , - t \cos \theta + s \sin \theta ) \, \rmd t \label{eq:RadonDef}
\end{eqnarray}
In sparse angle tomography, only a small number of such angular projections can be measured. In our study, we chose only $m_\theta = 45$ angles, evenly distributed in $[0,\pi)$. In addition, for a given angle $\theta_i$, we practically only measure the integrals of $\Radon[\uInf](\theta_i, s)$ over small $s$-intervals representing an array of $m_s = 500$ equal sized sensor pixels. In total, this leads to $m$ = $m_s \cdot m_\theta= 18.000$ measurements. The forward operator $A$ corresponds to the exact discretization of this measurement with respect to the pixel basis: All the operations involved in the measurement can be computed explicitly for indicator functions of rectangular sets. Further details of this step can be found in Section 2.3 in \cite{Lu14}. \\
The unknown function $\uInf$ to recover is a slightly scaled version of the \emph{Shepp-Logan phantom} \cite{ShLo74}, a toy model of the human head defined by 10 ellipses. Figure \ref{subfig:PhanCTReal} shows $\uInf$ and Figure \ref{subfig:PhanCTSparseData} the measurement data generated by discretizing $\uInf$ with a $768 \times 768$ pixel grid. We will refer to this scenario as "Phantom-CT".

\subsection{Accuracy Assessment} \label{subsec:AccAss}

\begin{figure}[tb]
\centering
\subfloat[][]{\includegraphics[width = 0.5\textwidth]{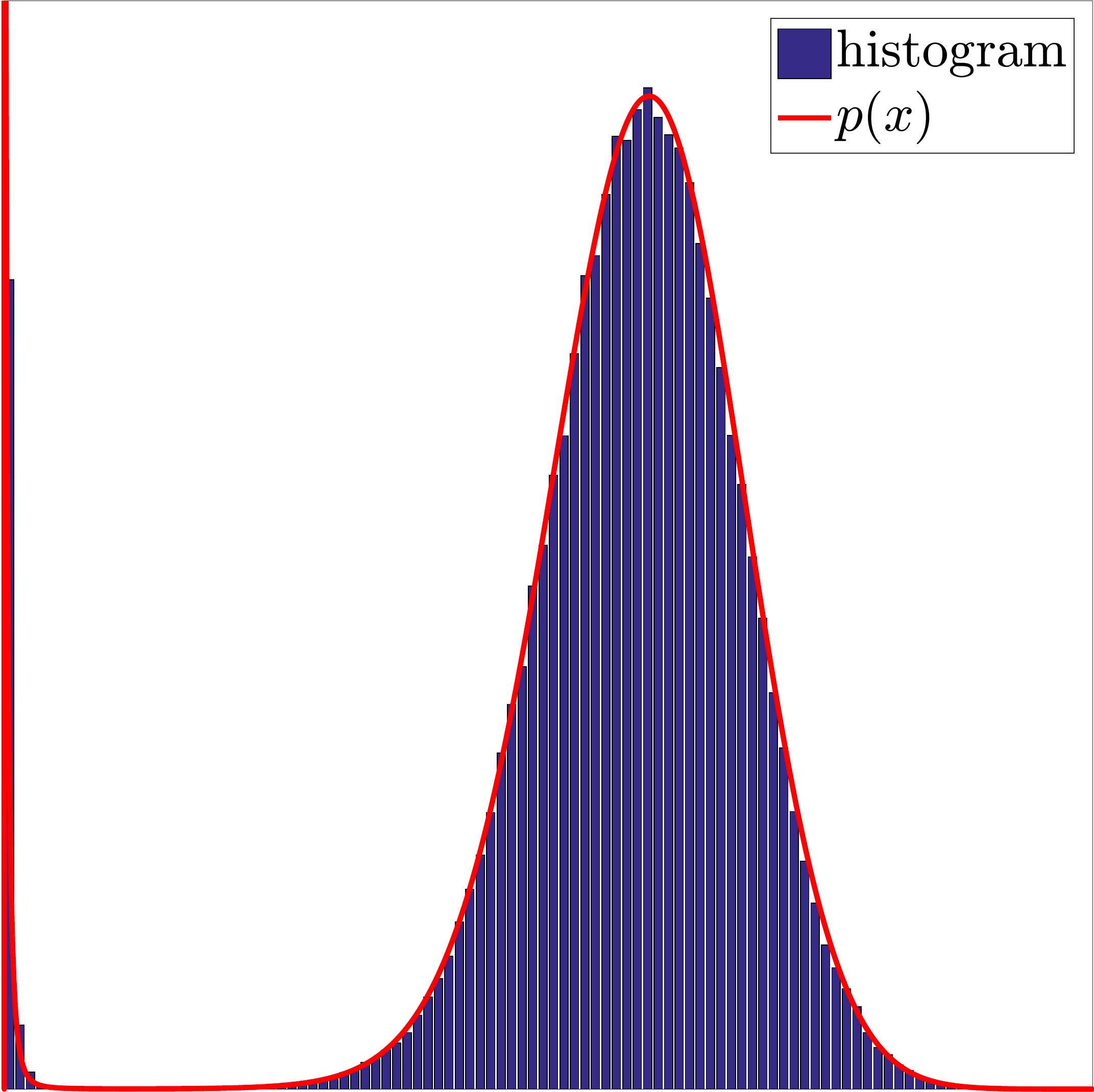} \label{subfig:SSevalImg1}}
\caption{\protect\subref{subfig:SSevalImg1} Histogram (blue bars) of the slice sampler compared to targeted SC density (red line) given by \eref{eq:SCLp} with $p = 0.8$. The parameters $a$ and $b$ were picked from a run of the direct $\ell_1$ sampler applied to a 2D image deblurring scenario (described in \cite{Lu12}) and $c$ was matched to the regularization parameter used therein. \label{fig:SSeval}
}
\end{figure}

To validate that the developed slice sampler accurately reproduces the distributions it is supposed to sample  from, the convergence of the sample histograms to the underlying SC densities was checked for visually. Various (random) combinations of coefficients for the different SC densities were tested; see Figure \ref{fig:SSeval} for an example of such a comparison. 

\subsection{Efficiency Assessment} \label{subsec:EffAss}

\begin{table}
\caption{\label{tbl:DirectVsSSTV} Comparison of $\tau_{int}$ for direct and slice-within-Gibbs samplers using different burn-in lengths for the slice sampler. The  "Boxcar" scenario and a TV prior ($p=q=1$) with $\lambda = 400$ is used and $K = 5\cdot 10^6$, $SSR = n$.} 
\lineup
\begin{tabular}{@{}*{6}{c}}        
\br                     
$K_{s,0} = 10$ & $K_{s,0} = 20$ & $K_{s,0} = 40$ & $K_{s,0} = 100$ & $K_{s,0} = 200$ & direct\\
\mr
231.4$\pm$8.6   & 149.2$\pm$4.6   & 109.4$\pm$2.9   & 102.0$\pm$2.6    & 101.3$\pm$2.6    & 97.8$\pm$2.5 \\
\br
\end{tabular}
\end{table}

\begin{table}
\caption{
\label{tbl:SS} Comparison of $\tau_{int}$ for slice-within-RSG samplers using different burn-in lengths for the slice sampler. The "Boxcar" scenario is used in (a) with an $\ell_p$ increment prior with $ p = 1.2$, $\lambda = 400$ and in (b) with an $\ell_p^q$ increment prior with $ p = 1$, $q = 10$, $\lambda = 0.02$. In both cases $K= 2 \cdot 10^6$ samples were drawn. In (c), the "Phantom-CT" scenario  with an isotropic TV prior with $\lambda = 500$ is used and $K = 2 \cdot 10^5$ samples were drawn.
} 
\begin{tabular}{@{}*{8}{c}}
\br                             
(a) \label{subtbl:DirectVsSSLp}  & $K_{s,0} = 1$ & $K_{s,0} = 2$ & $K_{s,0} = 4$ & $K_{s,0} = 8$ & $K_{s,0} = 16$ & $K_{s,0} = 32$ & $K_{s,0} = 64$ \\ 
\mr
  & 41.9$\pm$1.1 & 33.3$\pm$0.8 &23.4$\pm$0.5 & 18.3$\pm$0.3 & 15.8$\pm$0.4  & 14.6$\pm$0.3 & 14.8$\pm$0.3 \\
\br
\end{tabular}
\vskip 4pt
\begin{tabular}{@{}*{8}{c}}
\br                             
\label{subtbl:DirectVsSSLq} (b) & $K_{s,0} = 1\:$ & $K_{s,0} = 2\:$ & $K_{s,0} = 4\:$ & $K_{s,0} = 8\:$ & $K_{s,0} = 16\:$ & $K_{s,0} = 32\:$ & $K_{s,0} = 64$ \\ 
\mr
& 638$\pm$46 & 425$\pm$26 & 307$\pm$16 &  198$\pm$9 & 161$\pm$6  & 155$\pm$7 & 135$\pm$6 \\
\br
\end{tabular}
\vskip 4pt
\begin{tabular}{@{}*{8}{c}}
\br                             
\label{subtbl:SSevalTV} (c) & $K_{s,0} = 0$ & $K_{s,0} = 1$ & $K_{s,0} = 2$ & $K_{s,0} = 4$ & $K_{s,0} = 8$ & $K_{s,0} = 16$ & $K_{s,0} = 32$  \\ 
\mr
& 6.0$\pm$0.3 & 5.3$\pm$0.3 & 5.3$\pm$0.3 & 5.6$\pm$0.3 &  5.2$\pm$0.3 & 4.9$\pm$0.3  & 5.2$\pm$0.3 \\
\br
\end{tabular}
\end{table}

Once the accuracy of the slice sampler is established, the next crucial question is whether its use within a Gibbs sampler is \emph{efficient}: In Algorithm \ref{algo:Gibbs}, we ideally want to replace the current values of the component $j$,  $u_i$ by a values that is both distributed following $p_{post}(  \, \cdot \, | u_{-j},f)$ and statistically independent of the current value $u_i$. While direct SC samplers, such as the iCDF, naturally fulfill these requirements, SC samplers relying on MCMC chains initialized with $u_i$ fulfill them only asymptotically, in the limit $K_{s,0} \rightarrow \infty$. Using a fixed chain size $K_{s,0}$ will inevitably introduce additional correlation between subsequent samples and lower the \emph{statistical efficiency} of slice-within-Gibbs samplers compared to Gibbs sampling relying on a direct sampler for the SC densities. In the following, we will asses this loss of statistical efficiency by autocorrelation analysis. 

\paragraph*{Autocorrelation Analysis}
Evaluating samplers in general rather than for a specific aim is a difficult task \cite{Li08}. For the sake of a concise presentation, a detailed introduction and discussion is omitted here but can be found in Section 4.1.6. of \cite{Lu14}. In this study, we will only examine the autocorrelation functions of the MCMC chains projected onto a test function $v \in \R^n$, i.e., of the chain 
\begin{equation}
\{ g^i \}_{i=1}^{K} = \{ \langle v , u^i \rangle \}_{i=K_0+1}^{K_0+K}.
\end{equation}
In the "Boxcar" scenario, $v$ is given as the largest eigenvector of the (pre-computed) posterior covariance matrix while in the "Phantom-CT" scenario, it is the indicator function of the area defined by $[-0.32,-0.12] \times [0.12,0.32]$ (this area corresponds to the green box shown in Figures \ref{subfig:PhanTVMAP8}-\ref{subfig:PhanTVCM8}). To extract a quantitative measure from the autocorrelation functions, we will estimate their \emph{integrated autocorrelation time} $\tau_{int}$ by the approach presented in \cite{Wo04}. In all computed examples, the component $j$ to update in step A1.1 of Algorithm \ref{algo:Gibbs} is drawn uniformly at random, (\emph{random scan Gibbs sampler}) and a \termabb{sub-sampling rate}{SSR} of $n$ is used, i.e., only every $n$-th sample of the chains generated by Algorithm \ref{algo:Gibbs} is actually stored and $\tau_{int}$ refers to the samples of this thinned chain. This means that, on average, we update all $n$ components of $u \in \R^n$ between two steps of the chain (\emph{full sweep}). In each scenario, the samplers were given a large burn-in time $K_0$ and $K$ was chosen large enough to obtain sufficiently tight error bounds on $\tau_{int}$ \cite{Wo04}.

\paragraph*{Results}

When using a conventional TV prior ($p=q=1$) in the "Boxcar" scenario, the direct $\ell_1$ sampler using the iCDF method can be used as a reference to which the slice samplers can be compared to: The $\tau_{int}$ obtained by the direct sampler is a lower bound for all slice samplers. Table \ref{tbl:DirectVsSSTV} lists the results. One can observe that already for small MCMC chain length $K_{s,0}$, the differences between direct and slice sampler in terms of statistical efficiency are negligible in practice. Similar examinations using $\ell_2$ priors (where, again, a direct sampler can be used as a reference) showed that in this case, significant differences vanish for even smaller values of $K_{s,0}$ (results omitted here).\\
Tables \ref{tbl:SS} (a), (b) and (c) show the results of similar examinations for an $\ell_p$ prior with $p = 1.2$, an $\ell_p^q$ prior with $p = 1$, $q = 10$ and the isotropic TV prior in the 2D "Phantom-CT" scenario (using $n = 129 \times 129$), respectively. While we do not have a direct sampler as a reference here, one can clearly see that $\tau_{int}$ is converging to a limit for increasing $K_{s,0}$. In some cases, even using $K_{s,0} = 0$, i.e., only performing one step of the slice sampler, might be sufficient. 

\paragraph*{Computational Complexity}

In typical large scale inverse problems such as the one examined in Section \ref{subsec:CTTV}, the computational bottleneck is to compute the coefficients of the SC densities, not the process of sampling from them. Therefore, the computational complexity of the slice sampler is not a critical aspect of the whole algorithm. However, to give an indication of how increasing $K_{s,0}$ effects the total run time, Table \ref{tbl:DirectVsSSTVcomp} compares the run time of the slice-within-Gibbs sampler to the direct $\ell_1$ sampler and Table \ref{tbl:CompEffTV} lists the run times of the slice-within-Gibbs sampler for the TV prior in 2D. While the implementation of the slice sampling part is more complicated in this situation (cf. Section \ref{subsec:ImplTVSS}), we see that even for a moderate sized scenario ($n = 256^2$) it does not significantly effect the run time. Therefore, one does not have to compromise statistical efficiency by choosing a small $K_{s,0}$ to obtain a better computational efficiency. 

\begin{table}
\caption{\label{tbl:DirectVsSSTVcomp} Total run time of the slice-within-Gibbs sampler using different burn-in lengths divided by the run time of the direct $\ell_1$ sampler. The "Boxcar" scenario and an $\ell_p$ increment prior with $ p = 1.2$, $\lambda = 400$ is used.} 
\lineup
\begin{tabular}{@{}*{5}{c}}        
\br                     
$K_{s,0} = 10$ & $K_{s,0} = 20$ & $K_{s,0} = 40$ & $K_{s,0} = 100$ & $K_{s,0} = 200$ \\
\mr
1.30   & 1.38   & 1.48   & 1.75    & 2.20     \\
\br
 \end{tabular}
\end{table}

\begin{table}
\caption{\label{tbl:CompEffTV} Total run time of the slice-within-Gibbs sampler using different burn-in lengths $K_{s,0}$ divided by the run time for $K_{s,0} = 0$. The "Phantom-CT" scenario ($n = 256\times256$) and a TV prior ($p=1$) with $\lambda = 500$ is used.} 
\lineup
\begin{tabular}{@{}*{7}{c}}        
\br                     
$K_{s,0} = 0$  & $K_{s,0} = 1$ & $K_{s,0} = 2$ & $K_{s,0} = 4$ & $K_{s,0} = 8$ & $K_{s,0} = 16$ & $K_{s,0} = 32$ \\
\mr
1.00   & 0.98   & 1.04   & 1.03    & 1.05  &  1.07  &  1.11  \\
\br
 \end{tabular}
\end{table}

\subsection{Application to the  $\ell_p$ Increment Prior} \label{subsec:TVp}

\begin{figure}[tb]
\centering
\subfloat[][MAP estimates by SA]{\includegraphics[width = 0.48\textwidth]{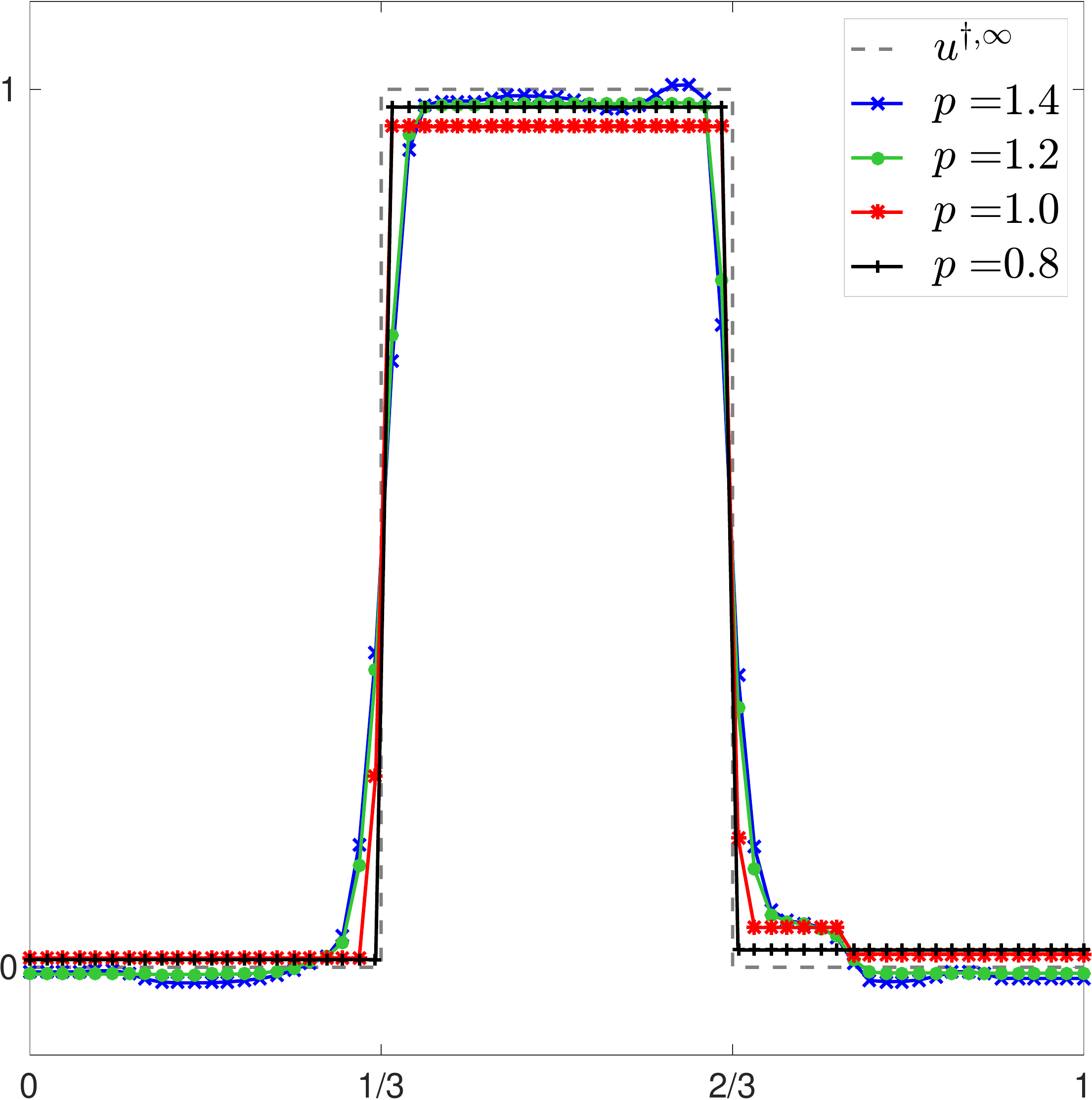} \label{subfig:BoxcarTVpMAP}}
\hspace{0.01\textwidth}
\subfloat[][CM estimates by RSG]{\includegraphics[width = 0.48\textwidth]{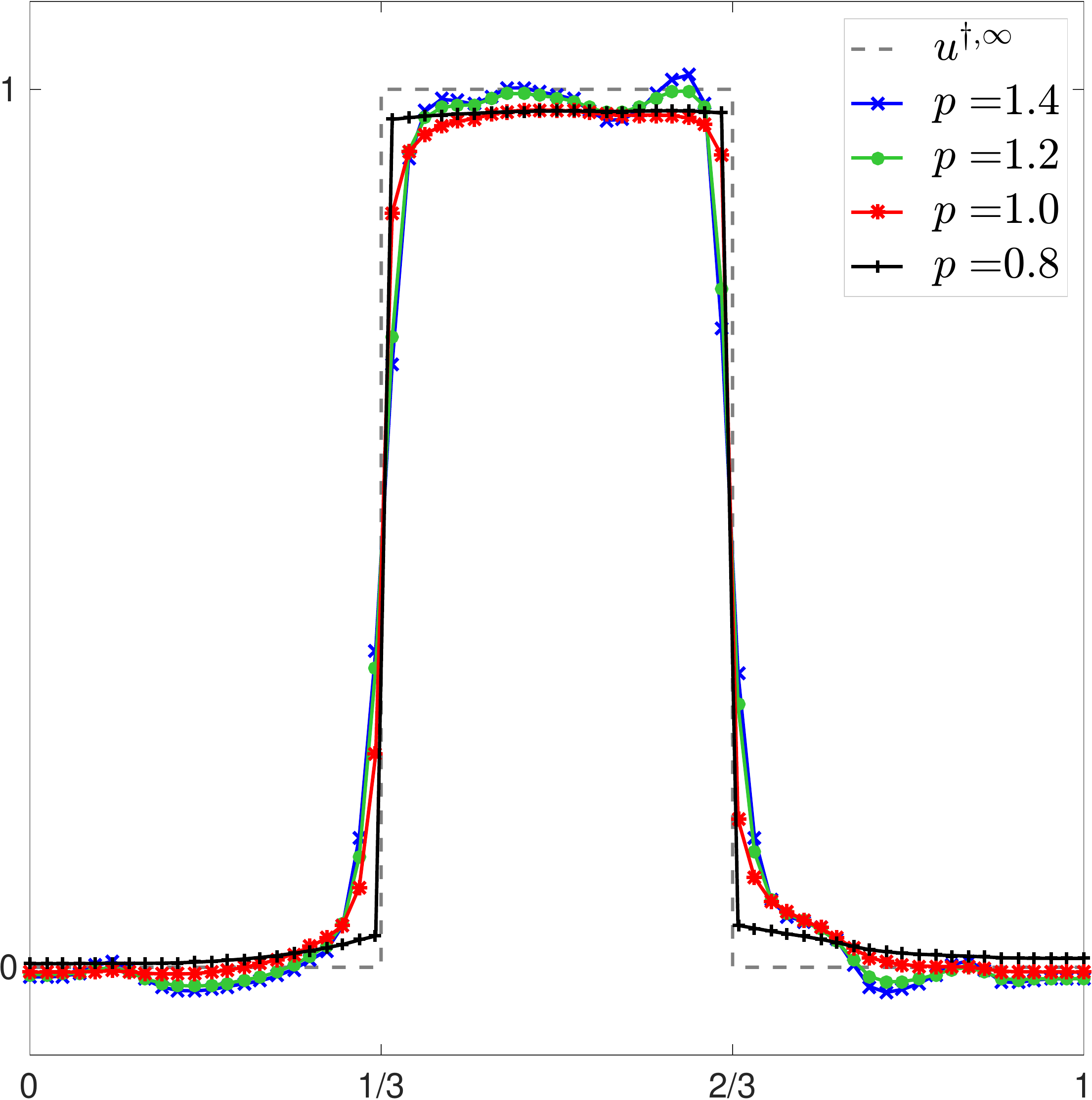} \label{subfig:BoxcarTVpCM}}
\caption{MAP and CM estimates for the 1D "Boxcar" scenario using the $\ell_p$ increment prior and $n = 63$. The true function to recover (gray line plot) is denoted by $u^{\dagger,\infty}$. \label{fig:BoxcarTVp}}
\end{figure}

Problems with the conventional TV prior (see \cite{LaSi04} and the overview in Section 5.2.2 in \cite{Lu14}) stimulated research into alternative, edge-preserving prior models. Here, we exemplify how the new slice-within-Gibbs sampler can be used to investigate such general questions in Bayesian inversion: We use it to compute both MAP and CM estimates for the $\ell_p$ increment prior with $p$ decreasing from $p=2$ (Gaussian prior) to $p=1$ (TV prior) and even below $p < 1$ (non-logconcave prior). While the computation of the CM estimates is straight-forward, computing MAP estimates is done by using the sampler within a \termabb{simulated annealing}{SA} scheme, a stochastic meta-heuristic for global optimization. The details and an evaluation of using SA together with the proposed SC Gibbs samplers can be found in Section 4.2.4 and 5.1.5 in \cite{Lu14}. In both cases, $K = 10^5$ samples were drawn with SSR = $n$.\\
The results of computing MAP and CM estimates for different values of $p$ are shown in Figure \ref{fig:BoxcarTVp}. Here, $\lambda$ was chosen such that all likelihood energies are equal and that $\lambda = 200$ for $p = 1$. The results suggest that using $p < 1$ leads to superior results for both MAP and CM estimates compared to $p = 1$. The MAP estimate is closer to the real solution as it is both sparser in the increment basis and the contrast loss is reduced. The CM estimate for $p = 0.8$ looks way more convincing compared to those for $p \geqslant 1$: It has clear pronounced edges that separate smooth, denoised parts. However, using the slice-within Gibbs samplers for $p < 1$ needs to be examined more carefully: While the results are visually convincing, we cannot be sure that the sampler explored the whole, possibly multimodal posterior and did not get stuck in a single mode.

\subsection{Application to CT Inversion with TV Priors} \label{subsec:CTTV}

\begin{figure}[tb]
\centering
\subfloat[][MAP, $n = 64 \times 64$]{\includegraphics[width = 0.45\textwidth]{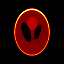} \label{subfig:PhanTVMAP6}}
\hspace{0.01\textwidth}
\subfloat[][CM, $n = 64 \times 64$]{\includegraphics[width = 0.45\textwidth]{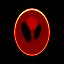} \label{subfig:PhanTVCM6}}\\
\subfloat[][MAP, $n = 128 \times 128$]{\includegraphics[width = 0.45\textwidth]{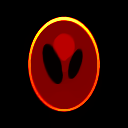} \label{subfig:PhanTVMAP7}}
\hspace{0.01\textwidth}
\subfloat[][CM, $n = 128 \times 128$]{\includegraphics[width = 0.45\textwidth]{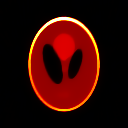} \label{subfig:PhanTVCM7}}\\
\subfloat[][MAP, $n = 256 \times 256$]{\includegraphics[width = 0.45\textwidth]{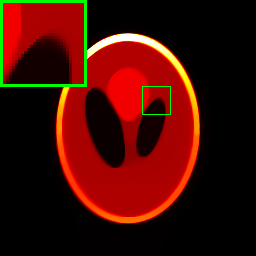} \label{subfig:PhanTVMAP8}}
\hspace{0.01\textwidth}
\subfloat[][CM, $n = 256 \times 256$]{\includegraphics[width = 0.45\textwidth]{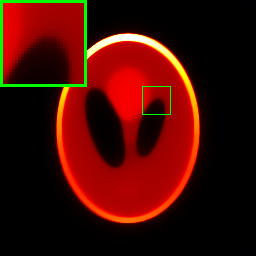} \label{subfig:PhanTVCM8}}\\
\caption{MAP and CM estimates for the "Phantom-CT" scenario using an isotropic TV prior with $\lambda = 500$, computed with increasing spatial resolution. In the highest resolution, a zoom inset is added. \label{fig:PhanTVMAPvsCM}}
\end{figure}

Figure \ref{fig:PhanTVMAPvsCM} shows MAP and CM estimates for the "Phantom-CT" scenario using an isotropic TV prior with Neumann boundary conditions \eref{eq:iTVprior}. Here, the MAP estimates were computed with the \termabb{alternating direction method of multipliers}{ADMM} \cite{BoPaChPeEc11}. In the Gibbs sampler, \termabb{oriented over-relaxation}{OOR} (see \cite{Ne95} and Section 4.3.1. in \cite{Lu14}) was used to accelerate convergence and $K = 2.5 \cdot 10^4, 10^4, 5 \cdot 10^3$ samples were drawn for $n = 64^2, 128^2, 256^2$, respectively (SSR = $n$).  

\subsubsection*{Non-Negativity Constraints}

\begin{figure}[tb]
\centering
\subfloat[][CM]{\includegraphics[width = 0.45\textwidth]{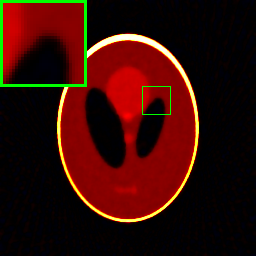} \label{subfig:PhanTVCM50}}
\subfloat[][CM, non-negativity const.]{\includegraphics[width = 0.45\textwidth]{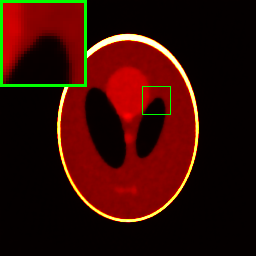} \label{subfig:PhanTVCM50nn}}
\\
\subfloat[][CStd]{\includegraphics[width = 0.45\textwidth]{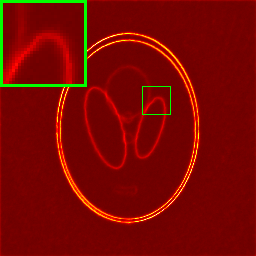} \label{subfig:PhanTVCStd50}}
\subfloat[][CStd, non-negativity const.]{\includegraphics[width = 0.45\textwidth]{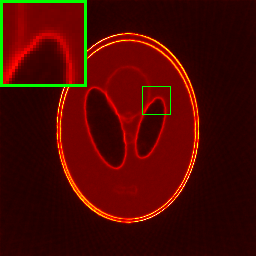} \label{subfig:PhanTVCStd50nn}}
\caption{Influence of inculding non-negativity constraints on CM and CStd estimates in the "Phantom-CT" scenario using an isotropic TV prior with $\lambda = 50$ and a spatial resolution of $n = 256 \times 256$. A zoom inset is added and both CM and both CStd estimates share the same color scale. \label{fig:PhanTV2}}
\end{figure}

As the slice-within-Gibbs sampler can easily incorporate additional hard constraints \eref{eq:HardCon}, it can be used to quantify their effect on the posterior $\post$. Figure \ref{fig:PhanTV2} shows CM and CStd estimates computed with or without  non-negativity constraints, $u \geqslant 0$. While the CM estimates look very similar, the CStd estimates reveal that the non-negativity constraints lead to a significant reduction of the posterior variance in some regions.

\subsubsection*{Gradient Estimates}

\begin{figure}[tb]
\centering
\subfloat[][CM]{\includegraphics[width = 0.45\textwidth]{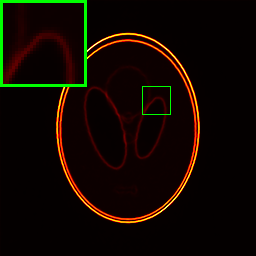} \label{subfig:PhanTVGrCM8}}
\subfloat[][CStd]{\includegraphics[width = 0.45\textwidth]{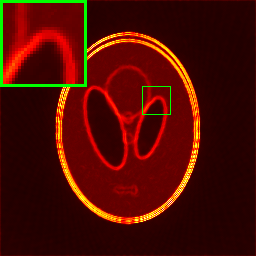} \label{subfig:PhanTVGrCStd8}}
\caption{CM and CStd estimates of $\norm{\nabla u}_2$ for the "Phantom-CT" scenario using an isotropic TV prior with $\lambda = 50$, non-negativity constraints and a spatial resolution of $n = 256 \times 256$. A zoom inset is added. \label{fig:PhanTV3}}
\end{figure}

We further present one example of how the samples of $u$ generated by the sampler can be used to compute statistics and uncertainties of a feature of $u$ (cf. Section\ref{subsec:SampleBasedInference}): In Figure \ref{subfig:PhanTVGrCM8}, we computed the CM estimate of the gradient of the image $u$,
\begin{equation}
\Exp \left[ \norm{ \nabla u}_2 |f \right] = \int \norm{ \nabla u}_2 \: \post \, \rmd u,
\end{equation}
and Figure \ref{subfig:PhanTVGrCStd8} shows the corresponding CStd estimate.


\section{Discussion and Conclusions}\label{sec:Discussion}

In this article, we presented and evaluated a new MCMC sampler that allows to carry out sample-based Bayesian inversion for a wide range of scenarios and prior models. It is based on the extension of the single component Gibbs-type sampler developed in \cite{Lu12} by a problem-specific adaptation and implementation of  generalized slice sampling and enables efficient posterior sampling in high-dimensional scenarios with certain priors for the first time.\\
The results in Sections \ref{subsec:AccAss} and \ref{subsec:EffAss} show that using generalized slice sampling to sample from the one-dimensional conditional, single component densities can lead to a fast, robust and accurate posterior sampler for the inverse problems scenario \eref{eq:FwdEq} and is therefore an attractive option whenever a fast direct sampler such as iCDF is not available. The computed results in Section \ref{subsec:TVp} exemplified the use of the new slice-within-Gibbs sampler to examine recent topics in Bayesian inversion and Section \ref{subsec:CTTV} demonstrated how it can lead to interesting results for Bayesian estimation in challenging, high-dimensional inverse problems scenarios. In particular, we examined that TV prior model in 2D: The theoretical analysis of the TV prior carried out, e.g. in \cite{LaSi04,LaSaSi09}, is restricted to 1D, only, and, to the best of our knowledge, no theoretical results are available for higher dimensions, yet. The development of the slice-within-Gibbs sampler now enabled us to examine the use the TV prior for the important inverse problem scenario of CT inversion in 2D, for the first time. The results show that, contrary to the 1D case, the CM estimates seem to get smoother for a constant value of $\lambda$ as the resolution increases. This observation could be the starting point for a new theoretical analysis and has to be examined in higher spatial dimensions by computational studies. \\
More generally, while our results and those of others (cf. Section \ref{subsec:SampleBasedInference}) have demonstrated that sampling high-dimensional posterior distributions is feasible for many important inverse problems scenarios nowadays, an important future challenge lies in extracting the information of interest from the samples generated: While we demonstrated, e.g., how to compute CStd estimates to examine how the posterior variance is influenced by non-negativity constraints (cf. Figure \ref{fig:PhanTV2}) or estimates of a feature $g(u)$ of $u$ (cf. Figure \ref{fig:PhanTV3}), we did not discuss how to interpret the corresponding results. This requires a concrete application and objective and will be topic of future investigations based on the methods presented here.\\
Related to the last point, we only used simulated data scenarios in this study to focus on the sampling algorithm. The application to experimental CT data will be the subject of a forthcoming publication covering more general aspects of Bayesian inversion in practical applications (see Section 5.3. in \cite{Lu14}). Furthermore, only prior models based on $\ell_p^q$-norms were considered here, while the sampler can, in principle, be implemented for more general prior models. A more fundamental limitation and future challenge is the current restriction to linear forward maps \eref{eq:FwdEq} and Gaussian noise models \eref{eq:Likelihood}. Both non-linear forward maps and non-Gaussian noise models typically conflict with condition (SC1), i.e., they make it very difficult to find an explicit parameterization of the SC densities. In addition, problems related to using SC-Gibbs sampling for multimodal posteriors (cf. Section \ref{subsec:TVp}) may occur as well.\\
Code to reproduce all the computed examples will be provided as part of the release of a Matlab-based toolbox for Bayesian inversion. 


\ack
This work was conducted as part of the PhD thesis \cite{Lu14} at the University of M\"unster. The author would like to thank the two supervisors, Martin Burger and Carsten H. Wolters. 

\appendix

\section{Details of the Implementation} \label{sec:Impl}

\subsection{Computation of the Likelihood Coefficients} \label{subsec:ImplLike}

To implement the SC Gibbs sampler for the "Boxcar" scenario, we compute $\Psi = A V$ and pre-compute $a \mydef \frac{1}{2} \sqnorm{\Psi_i}$. For computing $b = \Psi_i^T \varphi(i) = \Psi_i^T f - (\Psi_i^T \Psi_{-i}) \xi_{-i}$, we pre-compute $\Psi_i^T f$ and $\sqnorm{\psi_i }$ for all $i$ and build the $n \times n$ matrix $\Phi := \Psi^t  \Psi$. Then, computing $(\Psi_i^T \Psi_{-i}) \xi_{-i}$ can be performed by using 
\begin{equation}
 (\Psi_i^t \, \Psi_{[-i]}) \, \xi_{[-i]} = \xi^t \Phi_{(\cdot,i)} - \xi_i  \|  \psi_i \|^{2}_{2},
\end{equation}
which involves a scalar product of dimension $n$ as the most expensive operation.
\\
For the "Phantom-CT" scenario, there are two possible implementations: For the image sizes considered here ($n$ up to $256 \times 256$) we can still compute the matrices $\Psi = A$ and $\Phi := \Psi^t  \Psi$ explicitly and use the same implementation as in the "Boxcar scenario ($\Psi = A$ as we stay in the pixel basis, i.e., $V = I_n$). For larger $n$ or 3D applications, we might not be able to store $\Phi$ or $\Psi$. An alternative implementation that does not require storing any matrices uses
\begin{equation}
b = \Psi_{i}^T f - \Psi_{i}^T (\Psi \xi) + \xi_{i} \sqnorm{\Psi_{i}}
\end{equation} 
 to compute $b$ in the following way:
 \begin{itemize}
 \item We again pre-compute $\Psi_{i}^T f$ and $\sqnorm{\Psi_{i}}$ for all $i$. Then, we store the measurement that the current state $\xi$ would cause as $f_\xi$ and initialize it by $A u^0$. In principle, $f_\xi$ is given as $\Psi \xi$, and can be directly computed at any time but this computation is too expensive to be performed at every SC update. 
 \item For a given pixel $i$  that is to be updated, we construct $\Psi_{i}$ and compute the scalar product $\Psi_{i}^T f_\xi$ to update $b$ by the above formula  (note that $\Psi_{i}^T (f - \Psi \xi)$ is just a projection of $\Psi_{i}$ onto the current residual of $f_\xi = \Psi \xi$).  With the constructed $\Psi_{i}$ and the change, $\delta_{i}$, in $\xi_{i}$ caused by the sampling step,  we can then update $f_\xi = f_\xi + \delta_{i} \Psi_{i}$.
 \item While this iterative updating of $f_\xi$ is fast, inaccuracies can accumulate over time, leading to a misfit between $f_\xi$ and $\Psi \xi$. Therefore, we compute $\Psi \xi$ explicitly every $n$ steps and reset $f_\xi $ to this exact value.
 \end{itemize}
The computational bottleneck of this procedure is to compute $\Psi_{i}$, i.e., the Radon transform of a pixel (or voxel in 3D). For the parallel beam geometry used here, explicit formulas relying on basic operations that can be parallelized over the angles $\theta$ can be derived. For more complicated beam geometries, e.g., the
cone beam geometry for 3D reconstruction, approximations relying on basic operations from computer graphics can be derived and implemented very efficiently and parallelized on GPUs.

\subsection{Slice Sampling with TV Priors} \label{subsec:ImplTVSS}

\begin{figure}[tb]
\centering
\subfloat[][]{\includegraphics[width = 0.48\textwidth, height = 0.48\textwidth]{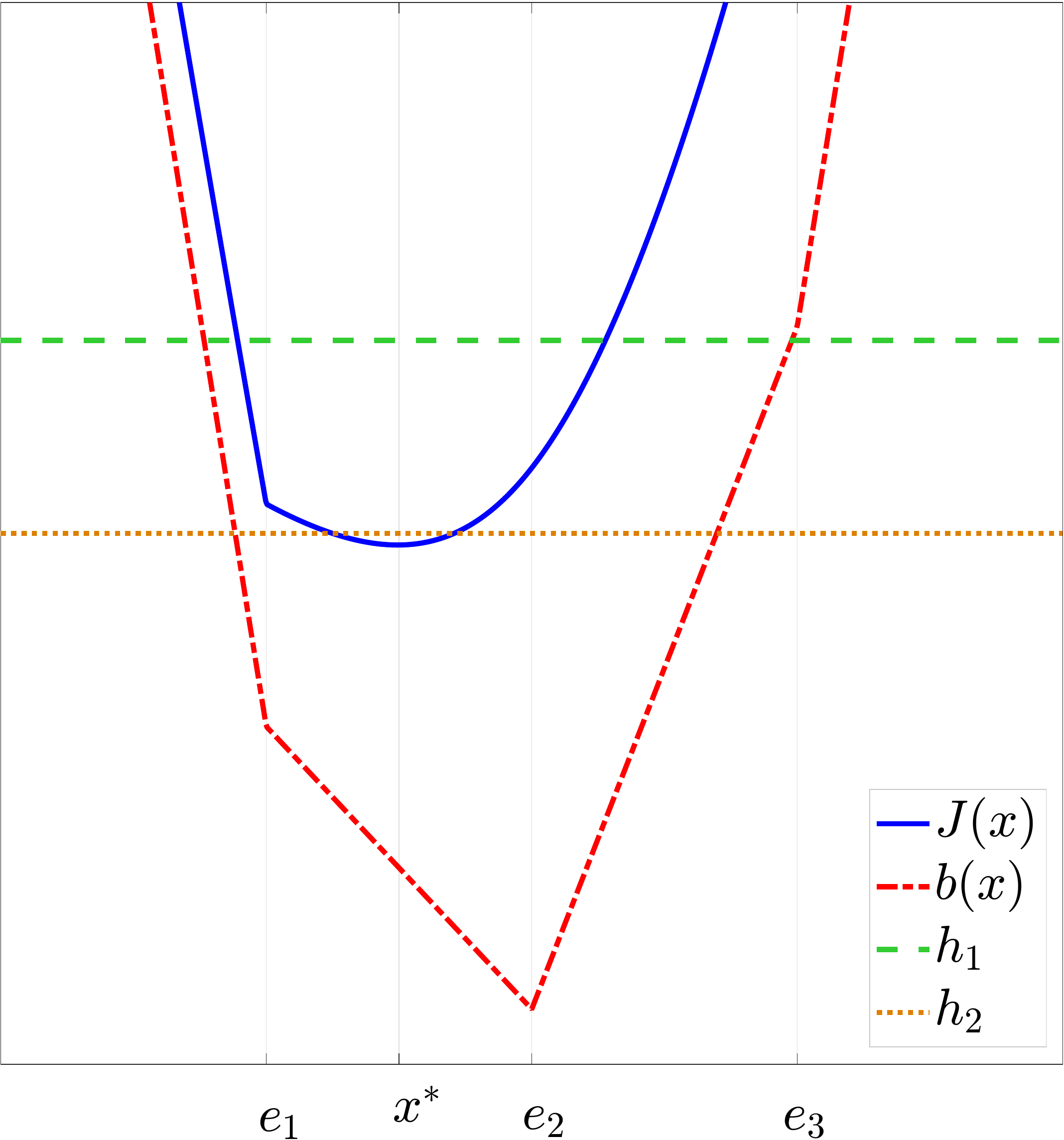}\label{subfig:TVSSImpl1}}
\hspace{0.01\textwidth}
\subfloat[][]{\includegraphics[width = 0.48\textwidth, height = 0.48\textwidth]{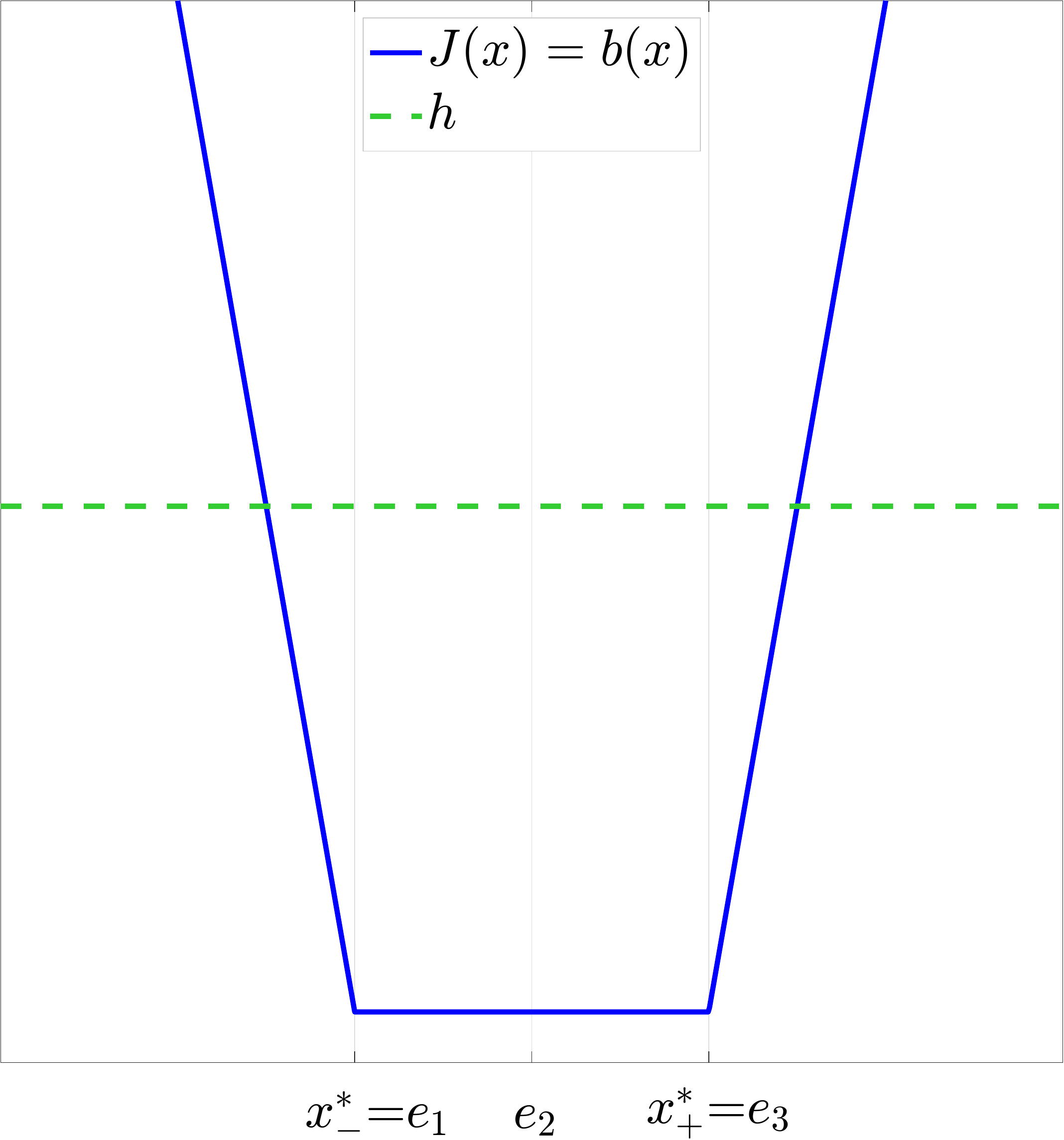} \label{subfig:TVSSImpl2}}
\caption[]{Illustration of two SC density energies for the slice sampler implementation of the TV prior in 2D: $J(x)$ (blue line), $b(x)$ (red line) and $h$ (green and yellow lines) for $(e_1,e_2,e_3) = (-1,0,1)$ and \subref{subfig:TVSSImpl1} $(d_1,d_2,d_3) = (2,1,1)$, $(g_1,g_2,g_3) = (0,0.5,1)$, \subref{subfig:TVSSImpl2} $(d_1,d_2,d_3) = (1,0,1)$, $(g_1,g_2,g_3) = (0,0,0)$\label{fig:TVSSImpl}.}
\end{figure}

From \eref{eq:SCTV}, we have 
\begin{equation}
\fl \quad p_2(x) = \exp \left( - c \sum_{j=1}^3 \sqrt{ d_j (x - e_j)^2  + g_j }\right), \quad d_j \in \{0,1,2\}, \quad  g_j \geqslant 0, 
\end{equation}
and want to determine $S_2^y = \left\lbrace z \;|\; p_2(z) \geqslant y \right\rbrace$ by solving   
\begin{equation}
y = p_2(x) \quad  \Longleftrightarrow \quad  - \frac{\log(y)}{c} =  \sum_{j=1}^3 \sqrt{ d_j (x - e_j)^2  + g_j } , 
\end{equation}
where $y \in (0,p_2(x))$ with probability 1 and $p_2(x) \leqslant 1$. Assume that $\left\lbrace e_1,e_2,e_3 \right\rbrace$ are sorted and define $J_j(x) \mydef \sqrt{ d_j (x - e_j)^2  + g_j }$ and $h \mydef - \log(y)/c$. Then, $J(x) \mydef \sum_j J_j(x)$ is convex and smooth in $I_1 \mydef (-\infty,e_1)$, $ I_2 \mydef (e_1,e_2)$, $ I_3 \mydef (e_2,e_3)$ and $I_4 \mydef  (e_3,\infty)$. It is monotonic in $I_1$ and $I_4$ and is bounded from below by  $b(x) \mydef \sum_j \sqrt{d_j} | x - e_j|$.  Define $\left[x_-^*,x_+^*\right] = \argmin \, J(x)$ as the interval of minimizers and $x_-$, $x_+$ as the solutions to $y = p_2(x)$. We have $x_- < x_-^*$,  $x_+ > x_+^*$, $x_- \in I_1 \cup I_2 \cup I_3$ and $x_+ \in I_2 \cup I_3 \cup I_4$ with probability 1 and $[x_-^*,x_+^*] \subset [e_1,e_3]$. See Figure \ref{fig:TVSSImpl} for two illustrations. \\
We will compute $x_-$ by a Newton's method:
\begin{equation}
x_-^i = x^{i-1}_{-} - \frac{J(x^{i-1}_{-}) - h}{J'(x^{i-1}_{-})},
\end{equation}
initialized in a point $x_-^0$ such that $x_-^0 \leqslant x_-$ and $J(x)$ is smooth on $[x_-^0,x_-]$. In each step, the Newton's method approximates $J(x)$ by a tangent in $x^{i-1}_-$. Due to the convexity of $J(x)$ and $x_-^0 \leqslant x_- < x_-^*$, the iterates never overshoot: $x_-^0 \leqslant x^i_- \leqslant x_-$ for all $i$. Thereby, they stay in $[x_-^0,x_-]$ and the derivative exists. Finding such an initialization $x_-^0$ requires some simple considerations: \\
The subdifferential $\partial J(x)$ is given as the sum of the subdifferentials of $J_i(x)$ (in the set-valued sense of addition):
\begin{equation}
\partial J_j(x) = \cases{
\left\lbrace \frac{d_j(x-e_j)}{\sqrt{ d_j (x - e_j)^2  + g_j }} \right\rbrace,  & if $x \neq e_j$ or $g_j > 0$\\
\left[-\sqrt{d},\sqrt{d} \right],  & if $x = e_j$ and $g_j = 0.$}
\end{equation}
Now, let $J^*_e \mydef \min_j J(e_j)$. We can distinguish two cases:
\begin{enumerate}[leftmargin=0.1\textwidth, itemsep=0pt,topsep=0pt]
\item[$h > J_e^*$:] In this case, we check the following conditions in sequence:
\begin{itemize}[leftmargin=0.03\textwidth , itemsep=0pt,topsep=0pt]
\item  If $h > J(e_1)$, $x_-$ is in $I_1$. We use the lower bound $b(x)$ to determine $x_-^0$ such that  $b(x_-^0) = h$:
\begin{equation}
	x_-^0 = e_1 +  \frac{J(e_1) - h}{\sum_j \sqrt{d_j}}
\end{equation}
 As $b(x) \leqslant J(x)$, and both are monotonic in $I_1$, we have that $x_-^0 < x_-$.
 \item Else if $h > J(e_2)$, $x_-$ is in $I_2$. We perform one Newton step from $e_1$ using the maximal subgradient in $e_1$:
   \begin{equation}
	x_-^0 = e_1 -  \frac{J(e_1) - h}{ \max \left( \partial J(e_1) \right) }
\end{equation}
This way, $x_-^0 \leqslant x_-$ and $[x_-^0,x_-] \subset I_2$, i.e., $J(x)$ is differentiable for all iterates.
\item Else, $h > J(e_3)$ and $x_-$ is in $I_3$. With a similar reasoning, we set 
   \begin{equation}
	x_-^0 = e_2 -  \frac{J(e_2) - h}{ \max \left( \partial J(e_2) \right) }
\end{equation}
\end{itemize}
For finding $x_+$, a similar reasoning can be applied. In the locations of non-differentiability, the minimal subgradient has to be used.
\item[$h < J_e^*$:] In this case, $J(x)$ is not piecewise linear (cf. the yellow line in Figure \ref{subfig:TVSSImpl1}) and the unique minimizer $x^*$ is not in  $\{e_1,e_2,e_3\}$. The convexity ensures that $x_- < x^* < x_+$ are all either in $I_2$ or $I_3$. If $\max(\partial J(e_1)) < 0$ and $\min(\partial J(e_2)) > 0$ we have that $x^*$ (and thereby $x_-$ and $x_+$) are in $I_2$. Otherwise, they are in $I_3$. As above, initial points $x_-^0$ and $x_+^0$ fulfilling the conditions can be found by performing one Newton step from the corners of the interval using the maximal/minimal subgradient. 
\end{enumerate}
The case $h = J^*_e$ has probability zero.  

\section*{References}

\bibliographystyle{abbrv}
\bibliography{all}

\begin{thebibliography}{10}

\bibitem{AgBaPaSt14}
S.~Agapiou, J.~M. Bardsley, O.~Papaspiliopoulos, and A.~M. Stuart.
\newblock {Analysis of the Gibbs Sampler for Hierarchical Inverse Problems}.
\newblock {\em SIAM/ASA Journal on Uncertainty Quantification}, 2(1):511--544,
  2014.

\bibitem{Ba12}
J.~M. Bardsley.
\newblock {MCMC-Based Image Reconstruction with Uncertainty Quantification}.
\newblock {\em SIAM Journal on Scientific Computing}, 34(3):A1316--A1332, 2012.

\bibitem{BaFo12}
J.~M. Bardsley and C.~Fox.
\newblock {An MCMC method for uncertainty quantification in nonnegativity
  constrained inverse problems}.
\newblock {\em Inverse Problems in Science and Engineering}, 20(4):477--498,
  2012.

\bibitem{BaSoHaLa14}
J.~M. Bardsley, A.~Solonen, H.~Haario, and M.~Laine.
\newblock {Randomize-Then-Optimize: A Method for Sampling from Posterior
  Distributions in Nonlinear Inverse Problems}.
\newblock {\em SIAM Journal on Scientific Computing}, 36(4):A1895--A1910, 2014.

\bibitem{BoPaChPeEc11}
S.~Boyd, N.~Parikh, E.~Chu, B.~Peleato, and J.~Eckstein.
\newblock {Distributed Optimization and Statistical Learning via the
  Alternating Direction Method of Multipliers}.
\newblock {\em Foundations and Trends in Machine Learning}, 3(1):1--122, Jan.
  2011.

\bibitem{BuLu14}
M.~Burger and F.~Lucka.
\newblock {Maximum a posteriori estimates in linear inverse problems with
  log-concave priors are proper Bayes estimators}.
\newblock {\em Inverse Problems}, 30(11):114004, 2014.

\bibitem{CaKaSo14}
D.~Calvetti, J.~P. Kaipio, and E.~Somersalo.
\newblock {Inverse problems in the Bayesian framework}.
\newblock {\em Inverse Problems}, 30(11):110301, 2014.

\bibitem{CaRoTa06}
E.~Candes, J.~Romberg, and T.~Tao.
\newblock {Robust uncertainty principles: exact signal reconstruction from
  highly incomplete frequency information}.
\newblock {\em IEEE Transactions on Information Theory}, 52(2):489--509, 2006.

\bibitem{Ch12}
N.~Chopin.
\newblock {Fast simulation of truncated Gaussian distributions}.
\newblock {\em Statistics and Computing}, 21(2):275--288, 2011.

\bibitem{Co11}
S.~Comelli.
\newblock {A Novel Class of Priors for Edge-Preserving Methods in Bayesian
  Inversion}.
\newblock Master's thesis, University of Milan, Italy, 2011.

\bibitem{CuFoSu11}
T.~Cui, C.~Fox, and M.~J. O'Sullivan.
\newblock {Bayesian calibration of a large-scale geothermal reservoir model by
  a new adaptive delayed acceptance Metropolis Hastings algorithm}.
\newblock {\em Water Resources Research}, 47, 2011.

\bibitem{CuLaMa16}
T.~Cui, K.~J. Law, and Y.~M. Marzouk.
\newblock Dimension-independent likelihood-informed mcmc.
\newblock {\em J. Comput. Phys.}, 304(C):109--137, Jan. 2016.

\bibitem{DaHaSt12}
M.~Dashti, S.~Harris, and A.~Stuart.
\newblock {Besov priors for Bayesian inverse problems}.
\newblock {\em Inverse Problems and Imaging}, 6:183--200, 2012.

\bibitem{Do06}
D.~L. Donoho.
\newblock {Compressed sensing}.
\newblock {\em IEEE Transactions on Information Theory}, 52(4):1289--1306,
  2006.

\bibitem{FoRa13}
S.~Foucart and H.~Rauhut.
\newblock {\em A Mathematical Introduction to Compressive Sensing}.
\newblock Birkh\"{a}user Basel, 2013.

\bibitem{HaLaMiSa06}
H.~Haario, M.~Laine, A.~Mira, and E.~Saksman.
\newblock {DRAM: Efficient adaptive MCMC}.
\newblock {\em Statistics and Computing}, 16(4):339--354, December 2006.

\bibitem{HaSaTa05}
H.~Haario, E.~Saksman, and J.~Tamminen.
\newblock {Componentwise adaptation for high dimensional {MCMC}}.
\newblock {\em Comput Stat}, 20(2):265--273, June 2005.

\bibitem{HmKaKoLaNiSi13}
K.~H\"am\"al\"ainen, A.~Kallonen, V.~Kolehmainen, M.~Lassas, K.~Niinimaki, and
  S.~Siltanen.
\newblock {Sparse Tomography}.
\newblock {\em SIAM Journal of Scientific Computing}, 35(3):B644--B665, 2013.

\bibitem{HeBu15}
T.~Helin and M.~Burger.
\newblock {Maximum a posteriori probability estimates in infinite-dimensional
  Bayesian inverse problems}.
\newblock {\em Inverse Problems}, 31(8):085009, 2015.

\bibitem{KaSo05}
J.~P. Kaipio and E.~Somersalo.
\newblock {\em {Statistical and Computational Inverse Problems}}, volume 160 of
  {\em {Applied Mathematical Sciences}}.
\newblock Springer New York, 2005.

\bibitem{KoLaNiSi12}
V.~Kolehmainen, M.~Lassas, K.~Niinim{\"a}ki, and S.~Siltanen.
\newblock {Sparsity-promoting Bayesian inversion}.
\newblock {\em Inverse Problems}, 28(2):025005 (28pp), 2012.

\bibitem{La14}
B.~Lanfer.
\newblock {\em Automatic Generation of Volume Conductor Models of the Human
  Head for EEG Source Analysis}.
\newblock PhD thesis, University of Muenster, 2014.

\bibitem{LaSaSi09}
M.~Lassas, E.~Saksman, and S.~Siltanen.
\newblock {Discretization invariant Bayesian inversion and Besov space priors.}
\newblock {\em Inverse Problems and Imaging}, 3(1):87--122, Feb 2009.

\bibitem{LaSi04}
M.~Lassas and S.~Siltanen.
\newblock {Can one use total variation prior for edge-preserving Bayesian
  inversion?}
\newblock {\em Inverse Problems}, 20:1537--1563, 2004.

\bibitem{LaRoRo13}
K.~Latuszynski, G.~O. Roberts, and J.~S. Rosenthal.
\newblock {Adaptive Gibbs samplers and related MCMC methods}.
\newblock {\em The Annals of Applied Probability}, 23(1):66--98, 02 2013.

\bibitem{Li08}
J.~Liu.
\newblock {\em {Monte Carlo Strategies in Scientific Computing}}.
\newblock {Springer Series in Statistics}. Springer New York, 2008.

\bibitem{Lu12}
F.~Lucka.
\newblock {Fast Markov chain Monte Carlo sampling for sparse Bayesian inference
  in high-dimensional inverse problems using L1-type priors}.
\newblock {\em Inverse Problems}, 28(12):125012, 2012.

\bibitem{Lu14}
F.~Lucka.
\newblock {\em Bayesian Inversion in Biomedical Imaging}.
\newblock PhD thesis, University of M\"unster, december 2014.

\bibitem{Ne95}
R.~M. Neal.
\newblock {Suppressing Random Walks in Markov Chain Monte Carlo Using Ordered
  Overrelaxation}.
\newblock Technical Report 9508, Learning in Graphical Models, 1995.

\bibitem{Ne03}
R.~M. Neal.
\newblock {Slice Sampling}.
\newblock {\em Annals of Statistics}, 31(3):705--767, 2003.

\bibitem{PaMaCh14}
A.~Palafox, M.~A. Capistran, and J.~A. Christen.
\newblock {Effective Parameter Dimension via Bayesian Model Selection in the
  Inverse Acoustic Scattering Problem}.
\newblock {\em Mathematical Problems in Engineering}, 2014:12, 2014.

\bibitem{Pe15}
M.~{Pereyra}.
\newblock {Proximal Markov chain Monte Carlo algorithms}.
\newblock {\em arXiv}, (1306.0187), 2015.

\bibitem{PeScChPeToHeMc15}
M.~{Pereyra}, P.~{Schniter}, E.~{Chouzenoux}, J.-C. {Pesquet}, J.-Y.
  {Tourneret}, A.~{Hero}, and S.~{McLaughlin}.
\newblock {A Survey of Stochastic Simulation and Optimization Methods in Signal
  Processing}.
\newblock {\em arXiv}, ({1505.00273}), 2015.

\bibitem{RoCa05}
C.~P. Robert and G.~Casella.
\newblock {\em Monte Carlo Statistical Methods (Springer Texts in Statistics)}.
\newblock Springer-Verlag New York, Inc., Secaucus, NJ, USA, 2005.

\bibitem{ShLo74}
L.~Shepp and B.~Logan.
\newblock The fourier reconstruction of a head section.
\newblock {\em IEEE Transactions on Nuclear Science}, 21(3):21--43, June 1974.

\bibitem{CoRoStWh13}
S.L.Cotter, G.O.Roberts, A.~Stuart, and D.~White.
\newblock Mcmc methods for functions: modifying old algorithms to make them
  faster.
\newblock {\em Statistical Science}, 28:424--446, 2013.

\bibitem{SoSo14}
S.~Sommariva and A.~Sorrentino.
\newblock {Sequential Monte Carlo samplers for semi-linear inverse problems and
  application to magnetoencephalography}.
\newblock {\em Inverse Problems}, 30(11):114020, 2014.

\bibitem{SoLuAr14}
A.~Sorrentino, G.~Luria, and R.~Aramini.
\newblock {Bayesian multi-dipole modelling of a single topography in MEG by
  adaptive sequential Monte Carlo samplers}.
\newblock {\em Inverse Problems}, 30(4):045010, 2014.

\bibitem{St10}
A.~M. Stuart.
\newblock {Inverse problems: A Bayesian perspective}.
\newblock {\em Acta Numerica}, 19:451--559, 5 2010.

\bibitem{Vo02}
C.~R. Vogel.
\newblock {\em Computational Methods for Inverse Problems}.
\newblock Society for Industrial and Applied Mathematics, Philadelphia, PA,
  USA, 2002.

\bibitem{Wo04}
U.~Wolff.
\newblock Monte carlo errors with less errors.
\newblock {\em Computer Physics Communications}, 156(2):143 -- 153, 2004.

\end{thebibliography}

\end{document}